# ON LARGE DEVIATION REGIMES FOR RANDOM MEDIA MODELS[1]

By M. Cranston, D. Gauthier and T. S. Mountford

*University of California, Irvine, École Polytechnique Fédérale de Lausanne and École Polytechnique Fédérale de Lausanne*

The focus of this article is on the different behavior of large deviations of random subadditive functionals above the mean versus large deviations below the mean in two random media models. We consider the point-to-point first passage percolation time $a_n$ on $\mathbb{Z}^d$ and a last passage percolation time $Z_n$. For these functionals, we have $\lim_{n\to\infty} \frac{a_n}{n} = \nu$ and $\lim_{n\to\infty} \frac{Z_n}{n} = \mu$. Typically, the large deviations for such functionals exhibits a strong asymmetry, large deviations above the limiting value are radically different from large deviations below this quantity. We develop robust techniques to quantify and explain the differences.

**1. Introduction.** We introduce the models to be treated in this work which is inspired by the results in [3, 5, 7, 9, 10, 11].

*Model 1.* The first model arises by analogy with the continuous time parabolic Anderson model. Define, for $x \in \mathbb{Z}^d, |x|_1 = |x_1| + |x_2| + \cdots + |x_d|$, and consider the graph $G = (\Xi, E)$ where

$$\Xi = \{(x,n) \in \mathbb{Z}^d \times \mathbb{Z}_+ : |x|_1 + n \equiv 0 \bmod 2\}$$

and $E$ denotes the set of directed nearest neighbor edges from vertices $(x, n)$ to vertices of the form $(x \pm e_i, n+1), i = 1, \ldots, d$, where $e_i = (0, \ldots, i, \ldots, 0)$ is the $i$th basis vector in $\mathbf{R}^d$. When $d = 1$, the edges connect vertices of the form $(x, n)$ to vertices of the form $(x \pm 1, n+1)$. Notice that $\underline{0} = (0,0) \in \Xi$ and the requirement $|x|_1 + n \equiv 0 \bmod 2$ implies there are $n+1$ elements of the form $(x, n) \in \Xi$ which are accessible from $\underline{0}$ and $2(n+1)$ upward directed edges from these points. As for notation, we shall use $x, y, z$ to denote elements of $\mathbb{Z}^d, \underline{x}, \underline{y}, \underline{z}$ to denote elements of $\Xi$ and $m, n, k, l$ to denote elements of $\mathbb{Z}_+$.

---

Received September 2006; revised December 2007.

[1]Supported in part by a NSF Grant DMS-07-06198 and SFNS 3510767.

Set

$$\Xi_n = \Xi \cap (\mathbb{Z}^d \times \{n\}) \tag{1}$$

and for $A \subset \mathbb{Z}^d$, set

$$\Xi_n(A) = \Xi \cap (A \times \{n\}). \tag{2}$$

We will call sets of the form $\Xi_n(A)$ blocks. If $I \subset \mathbb{Z}_+$, put

$$\Xi_I(A) = \bigcup_{n \in I} \Xi_n(A). \tag{3}$$

Define the set of nearest neighbor paths in $G$ from $\underline{x} = (x, k) \in \Xi_k$ to $\Xi_n$ by

$$\aleph_n^{\underline{x}} = \{\underline{\gamma} : [0, n-k] \to \Xi : \text{ so that } \underline{\gamma}(0) = \underline{x} \text{ (and so } \underline{\gamma}(n-k) \in \Xi_n)\}. \tag{4}$$

Since the edges in the graph $G = (\Xi, E)$ are directed, the $(d+1)$st coordinate of $\underline{\gamma}$ denoted $\underline{\gamma}_{d+1}$ will be increasing. We shall sometimes write $\underline{\gamma} = (\gamma, \underline{\gamma}_{d+1})$. If $e \in E$ is an edge along $\underline{\gamma}$, we shall write $e \in \underline{\gamma}$. We also introduce, for $m < n$,

$$\aleph_{m,n} = \bigcup_{\underline{x} \in \Xi_m} \aleph_n^{\underline{x}}. \tag{5}$$

For $A \subset \mathbb{Z}^d$, we set

$$\aleph_{m,n}^A = \bigcup_{\underline{x} \in \Xi_m(A)} \aleph_n^{\underline{x}}. \tag{6}$$

For simplicity, write $\aleph_n = \aleph_n^0$. Now consider an i.i.d. random field $\{X_e : e \in E\}$ defined on some probability space $(\Omega, \mathcal{F}, P)$ with $E[X_e] = 0$. One choice of interest for the distribution of $X_e$ is $\mathcal{N}(0,1)$. We shall, however, at the very least always impose the Cramér condition

$$E[e^{cX_e}] < \infty \qquad \forall |c| < c_0 \text{ for some } c_0 > 0. \tag{7}$$

We set

$$Z_n = \sup_{\underline{\gamma} \in \aleph_n} \sum_{e \in \underline{\gamma}} X_e. \tag{8}$$

Again, simple subadditive considerations lead us to the conclusion that there is a nonrandom constant $\mu$ such that

$$\lim_{n \to \infty} \frac{Z_n}{n} = \mu, \qquad P\text{-a.s.}$$

We will investigate the asymptotic behavior (for small $\varepsilon > 0$) of the upper and lower large deviation probabilities $P(Z_n > (\mu + \varepsilon)n)$ and $P(Z_n < (\mu - \varepsilon)n)$ as $n$ tends to infinity.

*Model 2.* The second model is standard first passage percolation and our results are motivated by the work of Chow and Zhang [2]. In that work, first passage percolation questions were considered. Here, we consider the graph $G = (\mathbb{Z}^d, E)$, where $E$ denotes the set of nearest neighbor edges in $\mathbb{Z}^d$. When $x$ and $y$ are adjacent edges in $\mathbb{Z}^d$, we shall denote the corresponding edge by $e = \{x, y\}$. Here, we consider the graph $G = (\mathbb{Z}^d, E)$ where $E$ denotes the set of nearest neighbor edges in $\mathbb{Z}^d$.

The set over which optimization will take place is the set, $\Theta$, of nearest neighbor paths $\gamma$ in $G$ having finite length denoted by $l(\gamma)$. Given $A, B \subset \mathbb{Z}^{d-1}, m < n$ write

$$(9) \quad \Theta_{m,n}(A, B) = \{\gamma \in \Theta : \gamma(0) \in A \times \{m\}, \gamma(l(\gamma)) \in B \times \{n\}\}.$$

For notational simplicity, we drop 0 in the case $m = 0$ writing

$$\Theta_n(A, B) = \Theta_{0,n}(A, B)$$

and further simplify in the case $A = [-n, n]^{d-1}, B = \mathbb{Z}^{d-1}$, by writing

$$\Theta_n = \Theta_n([-n, n]^{d-1}, \mathbb{Z}^{d-1})$$

for the paths which start in the box $[-n, n]^{d-1} \times \{0\}$ and terminate at the hyperplane $\mathbb{Z}^{d-1} \times \{n\}$. The random variable $G_n$ was defined in [2] as

$$(10) \quad G_n \equiv \inf_{\gamma \in \Theta_n} \sum_{e \in \gamma} t_e.$$

Similarly, if we again take nearest neighbor paths of finite length and define the subset,

$$(11) \quad \Psi_n = \Theta_n(\{0\}, \{(n, 0, \ldots, 0)\})$$
$$(12) \quad = \{\gamma \in \Theta : \gamma(0) = 0, \gamma(l(\gamma)) = (n, 0, \ldots, 0)\},$$

one sets

$$(13) \quad a_n = \inf_{\gamma \in \Psi_n} \sum_{e \in \gamma} t_e.$$

It was shown in [6] that there exists a finite nonrandom constant $\nu$, so that

$$\lim_{n \to \infty} \frac{G_n}{n} = \lim_{n \to \infty} \frac{a_n}{n} = \nu \geq 0, \qquad P\text{-a.s.}$$

Under the condition that the $t_e$ are bounded, it was shown in [8] that for some positive $A$ and $B$,

$$P(a_n \geq (\nu + \varepsilon)n) \leq A e^{-Bn^d}.$$

When the distribution of the $t_e$ satisfy $F(0) < p_c$, where $p_c$ is the critical percolation probability, then $\nu > 0$. Under the assumption $F(0) < p_c$, [2] showed that for small $\varepsilon > 0$, there are positive constants $c(\varepsilon), c'(\varepsilon)$ such that

$$\lim_{n\to\infty} \frac{-1}{n} \log P(G_n < (\nu - \varepsilon)n) = c(\varepsilon),$$

$$\lim_{n\to\infty} \frac{-1}{n^d} \log P(G_n > (\nu + \varepsilon)n) = c'(\varepsilon).$$

(Note: Throughout this paper, we use log to denote the logarithm base 2.)

The reasons for this quantitative disparity between upper and lower large deviation rates is that for $G_n$ to be small requires essentially one aberrant path $\gamma$ along which $\sum_{e \in \gamma} t_e$ is small, while for $G_n$ to be large requires all paths $\gamma$ to have $\sum_{e \in \gamma} t_e$ large. See [4] for similar phenomenon. This raises the question of the corresponding behavior for $P(\frac{a_n}{n} < \nu - \varepsilon)$ and $P(\frac{a_n}{n} > \nu + \varepsilon)$. As remarked in [2], the condition $E[e^{ct_e}] < \infty$ cannot ensure that (for $d > 1$)

$$\varlimsup_{n\to\infty} \frac{1}{n^d} \log P(a_n > (\nu + \varepsilon)n) < 0.$$

One seeks natural conditions on the distribution of $t_e$ so that

$$\varlimsup_{n\to\infty} \frac{1}{n^d} \log P(a_n > (\nu + \varepsilon)n) < 0.$$

However, provided $\varepsilon$ is small enough that $P(t_e > \nu + \varepsilon) > 0$, by considering the environment in which $t_e > \nu + \varepsilon$ for all $e \in [-n, n]^d$.

$$\varliminf_{n\to\infty} \frac{1}{n^d} \log P(a_n > (\nu + \varepsilon)n) > -\infty$$

is readily seen to hold.

Using the subadditivity arguments of [1, 8] we have that [provided that $P(t_e \leq \nu - \varepsilon) > 0$],

(14) $$-\infty < \lim_{n\to\infty} \frac{1}{n} \log P(a_n \leq (\nu - \varepsilon)n) < 0$$

and [provided that $P(X_1 \geq \mu + \varepsilon) > 0$]

(15) $$-\infty < \lim_{n\to\infty} \frac{1}{n} \log P(Z_n \geq (\mu + \varepsilon)n) < 0.$$

We first examine the influence of the tail of the distribution of the $X_e$ on the large deviation behavior of $Z_n$.

THEOREM 1.1.  *For Model 1, assume for some positive $M_0 < \infty$ that there is a positive increasing function $f$ so that*

$$\log P(X_e < x) = -(-x)^{d+1} f(-x), \qquad x < -M_0.$$

*Then for every sufficiently small $\varepsilon > 0$,*

$$\varlimsup_{n \to \infty} \frac{1}{n^{d+1}} \log P(Z_n \leq (\mu - \varepsilon)n) < 0$$

*if and only if*

(16) $$\sum_{n=1}^{\infty} \frac{1}{f(2^n)^{1/d}} < \infty.$$

We remark that an examination of the proof enables a weakening of the monotonicity assumption on $f$. It could be replaced by a condition such as the existence of finite constants $M_0$ and $c$ such that for $x < y < -M_0$, one has $cf(-x) > f(-y)$.

For Model 1, Theorem 1.1 implies that with $\{X_e : e \in E\}$ i.i.d. $\mathcal{N}(0,1)$ random variables, for some $\varepsilon > 0$

$$\lim_{n \to \infty} \frac{1}{n^{d+1}} \log P(Z_n \leq (\mu - \varepsilon)n) = 0$$

[and, in fact, for all $\varepsilon > 0$, $\varlimsup_{n \to \infty} \frac{1}{n^{d+1}} \log P(Z_n \leq (\mu - \varepsilon)n) = 0$]. However, we will now refine this result, giving the precise large deviation behavior for $\{X_e : e \in E\}$ i.i.d. $\mathcal{N}(0,1)$. In the Gaussian case, we have the following theorem.

THEOREM 1.2. *For Model 1 with the $\{X_e : e \in E\}$ i.i.d. $\mathcal{N}(0,1)$ and $d = 1$, for every sufficiently small $\varepsilon > 0$,*

$$-\infty < \varliminf_{n \to \infty} \frac{\log n}{n^2} \log P(Z_n \leq (\mu - \varepsilon)n) \leq \varlimsup_{n \to \infty} \frac{\log n}{n^2} \log P(Z_n \leq (\mu - \varepsilon)n) < 0.$$

If the above model is considered for $d > 1$ dimensions, then it can be seen that the lower large deviation have a rate of $n^2$, that is, there is no log term. Theorem 1.1 generalizes in an obvious way to the analogous question for Model 2 with $\{t_e : e \in E\}$ i.i.d. with $E$ the set of undirected nearest neighbor edges in $\mathbb{Z}^d$. Though this is a slightly different framework than the problem of Model 1, we see the following theorem.

THEOREM 1.3. *For Model 2, assume for some positive $M_0 < \infty$ that there is a positive increasing function $f$ so that*

$$\log P(t_e > x) = -x^d f(x), \qquad x > M_0.$$

*Then for every sufficiently small $\varepsilon > 0$,*

$$\varlimsup_{n \to \infty} \frac{1}{n^d} \log P(a_n \geq (\nu + \varepsilon)n) < 0$$

*if and only if*

$$\sum_{n=1}^{\infty} \frac{1}{f(2^n)^{1/(d-1)}} < \infty. \tag{17}$$

In this paper, we will be dealing with paths and the value of paths. The meaning of these words or phrases will depend on the model discussed.

In Model 1, a path $\underline{\gamma} = (\gamma, \gamma_{d+1})$ is a function from a finite interval, $I$, of nonnegative integers to $\Xi$, so that for all $n \in I, \gamma_{d+1}(n+1) = \gamma_{d+1}(n) + 1$ and $|\gamma(n+1) - \gamma(n)|_1 = 1$. If the interval is $I = [n_1, n_2]$, then the value of the path $\underline{\gamma}$ is

$$V(\underline{\gamma}) = \sum_{i=n_1}^{n_2-1} X_{e_i},$$

where the edge $e_i = (\underline{\gamma}(i), \underline{\gamma}(i+1))$. However, for simplicity, we shall use the notation

$$V(\underline{\gamma}) = \sum_{e \in \underline{\gamma}} X_e.$$

A path for Model 2 is simply a function $\gamma$ on a finite interval of positive integers, $I$, so that for all $n \in I, |\gamma(n) - \gamma(n+1)|_1 = 1$ (where both terms are defined). The value of $\gamma$ defined on $I = [n_1, n_2]$ is

$$V(\gamma) = \sum_{i=n_1}^{n_2-1} t_{e_i}$$

with the undirected edge $e_i = (\gamma(i), \gamma(i+1))$. Again, we simplify this by writing

$$V(\gamma) = \sum_{e \in \gamma} t_e.$$

In both cases, for a path $\gamma : [a, b] \to \mathbb{Z}^d$, $\gamma(a)$ is the initial or starting value and $\gamma(b)$ is the end value or endpoint.

Before embarking on the technical details, we make a basic remark on the idea behind the proof of the theorems. From now on, the discussion will be for the case $d = 1$, but the results extend easily to higher dimensions. We will see in Section 2, Proposition 2.1, that in the context of Model 1, under the assumption (7), for all $\delta \in (0, 1)$, if

$$H_{\delta n, n} = \{\underline{x} \in \Xi_{[\delta n]}(-\delta n/4, \delta n/4) : \exists \underline{\gamma} \in \aleph_n^{\underline{x}}, V(\underline{\gamma}) \geq (\mu - \varepsilon/10)(1-\delta)n\} \tag{18}$$

and

$$G_{\delta n, n} = \{|H_{\delta n, n}| \geq \tfrac{9}{10} |\Xi_{[\delta n]}(-\delta n/4, \delta n/4)|\}, \tag{19}$$

then for some $c > 0$ and all $n$ large enough

$$P(G^c_{\delta n, n}) \leq e^{-cn^2}.$$

That is, outside an event of probability $e^{-cn^2}$, starting from most points $\underline{x} \in \Xi_{[\delta n]}(-\delta n/4, \delta n/4)$ there are paths $\underline{\gamma}$ to the hyperplane $\Xi_n$ along which $V(\gamma)$ is approximating the supremum of the sum of $X_e$ over paths of this length. This gives control on the contribution far from the starting point. Define

(20)
$$\begin{aligned}J_{0,\delta n} = \{\underline{x} \in \Xi_{[\delta n]}(-\delta n/4, \delta n/4) : \exists \underline{\gamma} \in \aleph_{[\delta n]}, \\ \underline{\gamma}(\delta n) = \underline{x}, V(\gamma) \geq -\varepsilon n/10\}.\end{aligned}$$

If we can show, again under condition (7), that for

(21) $$K_{0,\delta n} = \{|J_{0,\delta n}| \geq \tfrac{9}{10}|\Xi_{[\delta n]}(-\delta n/4, \delta n/4)|\}$$

one has

$$P(K^c_{0,\delta n}) \leq e^{-cn^2},$$

then we would know that with probability at least $1 - 2e^{-cn^2}$ there would be a point $\underline{x} \in \Xi_{[\delta n]}(-\delta n/4, \delta n/4)$, a path $\underline{\gamma}_1 \in \aleph_{\delta n}$ with $\underline{\gamma}_1(\delta n) = \underline{x}$, a path $\underline{\gamma}_2 \in \aleph_n^{\underline{x}}$ such that both $V(\gamma_1) \geq -\varepsilon n/10$ and $V(\gamma_2) \geq (1-\delta)(\mu - \varepsilon/10)n$. This gives control on the contribution near the starting point. Concatenating $\underline{\gamma}_1$ and $\underline{\gamma}_2$ gives a path $\underline{\gamma} \in \aleph_n$ for which (if $\varepsilon > \frac{10\mu\delta}{10+\delta}$) we have $V(\gamma) \geq (\mu - \varepsilon)n$. Since $Z_n \geq V(\gamma)$, it would follow that for some $c$ not depending on $n$ (for $n$ sufficiently large)

(22) $$P(Z_n \leq (\mu - \varepsilon)n) \leq 2e^{-cn^2}.$$

The moral is that the most likely way for $\frac{Z_n}{n}$ to be small is for there to be "unavoidably" many large negative $X_e$'s around the point $(0,0)$. In the Gaussian case, the probability in (22) is larger due to the fatter tails of the $X_e$. The ideas for Model 2 are similar and are outlined in Section 5.

We end the Introduction by remarking that our arguments are sufficiently robust to extend to site percolation problems which are natural analogues.

**2. Analysis far from the starting point.** In this section, we will concretely discuss Model 1, but similar considerations apply to Model 2. Any differences will be written explicitly in the following. Our principal results will state that essentially (outside probability $e^{-cn^2}$) it is impossible to ensure that the functional $Z_n$ is small via aberrant values of $X_e$ for $e$ of order $n$ away from the point $(0,0)$. Thus, small values of $Z_n$ must arise from large negative deviations of the field $\{X_e : e \in E\}$ for $e$ comparatively close to the initial point. The points close to the initial point are "unavoidable" for paths whereas there are ample points far from the starting point which allow paths to avoid patches of large negative deviations of the $X_e$.

Recall $\Xi_n$ from (2) and the definitions of $H_{\delta n,n}, G_{\delta n,n}, J_{0,\delta n}$ and $K_{0,\delta n}$ from (18)–(20) and (21), respectively, of the previous section. We state and prove the following results in the case $d = 1$, but suitably stated analogs hold in all dimensions.

PROPOSITION 2.1. *Under (7), given any sufficiently small positive $\varepsilon$ and any $\delta \in (0,1)$, there is an $n_0$ and a constant $c = c(\varepsilon, \delta) > 0$ so that for all $n \geq n_0$ one has*

$$P(G^c_{\delta n,n}) \leq e^{-cn^2}.$$

REMARK. For general dimension $d$, with the analogously defined events, we will have the bound $P(G^c_{\delta n,n}) \leq e^{-cn^{d+1}}$.

The crucial part is a renormalization argument which follows using the reasoning of [1], Lemma 2.4.

DEFINITION 2.1. For a path $\underline{\gamma} \in \aleph_{m,n}$ and a subset $A \subset \mathbb{Z}^d$, we say $\underline{\gamma} \subset A$, if $\gamma(k) \in A, k = 0, 1, \ldots, n - m$. The set of such paths will be denoted by $\aleph_{m,n}(A)$. Similarly, if $\underline{x} \in A \times m$, then $\aleph_n^{\underline{x}}(A)$ denotes the set of paths $\underline{\gamma} \in \aleph_n^{\underline{x}}$ for which $\gamma(k) \in A, k = 0, 1, \ldots, n - m$.

LEMMA 2.1. *Given $\varepsilon > 0$ and $l \in \mathbb{Z}_+$ let*

$$A_l(\varepsilon) = \{\forall \underline{x} \in \Xi_0([0,l)), \exists \underline{\gamma} \in \aleph_l^{\underline{x}}([0,l)), V(\gamma) \geq (\mu - \varepsilon)l\}.$$

*There is a $c > 0$ such that given any $\varepsilon > 0$, there is an $l_0$ such that for $l \geq l_0$,*

$$P(A_l(\varepsilon)) \geq 1 - c\varepsilon.$$

PROOF. For any $\varepsilon > 0$, there exists $l_0$ so that for $l' \geq l_0$ (for notational convenience, we will suppose that $l'$ and $l'/\varepsilon$ are even integers)

$$P\left(\exists \underline{\gamma} \in \aleph_{l'}, V(\gamma) \geq \left(\mu - \frac{\varepsilon}{2}\right)l'\right) \geq 1 - \left(\frac{\varepsilon^3}{100}\right)^2.$$

From the FKG inequality applied to the decreasing events

$$A_1 = \left\{ \not\exists \underline{\gamma} \in \aleph_{l'}, \gamma(l') \in \mathbb{Z}_+, V(\gamma) \geq \left(\mu - \frac{\varepsilon}{2}\right)l' \right\},$$

$$A_2 = \left\{ \not\exists \underline{\gamma} \in \aleph_{l'}, \gamma(l') \in \mathbb{Z}_-, V(\gamma) \geq \left(\mu - \frac{\varepsilon}{2}\right)l' \right\},$$

we have

$$\left(\frac{\varepsilon^3}{100}\right)^2 \geq P(A_1 \cap A_2)$$
$$\geq P(A_1)P(A_2)$$
$$= P(A_1)^2$$

and consequently,

(23) $$P(A_1^c) = P(A_2^c) \geq 1 - \frac{\varepsilon^3}{100}.$$

Notice that the $\underline{\gamma}$ under consideration in either $A_1^c$ or $A_2^c$ must satisfy $\gamma \subset [-l', l']$ since such a $\underline{\gamma}$ starts at $\underline{0}$ and has only $l'$ steps of size 1. We now concatenate paths. First find, with high probability, a path $\underline{\gamma}_1 \in \aleph_{l'}$ for which $V(\gamma_1) \geq (\mu - \frac{\varepsilon}{2})l'$ by selecting a path as prescribed in $A_1^c$. This path satisfies both $\gamma_1 \subset [-l', l']$ and $\gamma_1(l') \geq 0$. Treat $(\gamma_1(l'), l')$ as the new origin, and continue by selecting a path $\underline{\gamma}_2 \in \aleph_{2l'}^{\gamma_1(l')}$ in an appropriately shifted version of $A_2^c$. By stationarity of the medium, the shifted versions of $A_1^c$ and $A_2^c$ also satisfy (23). Note this path stays in $[-2l', 2l']$ and has $\gamma_2(l') \in [-l', \gamma_1(l')] \subset [-l', l']$ and $V(\gamma_2) \geq (\mu - \frac{\varepsilon}{2})l'$. Repeating this procedure $\frac{2}{\varepsilon}$ times of going back and forth to stay in $[-2l', 2l']$ and concatenating the resulting paths gives a path $\underline{\gamma} \in \aleph_{(2l')/\varepsilon'}([-2l', 2l'])$ with $V(\gamma) \geq (\mu - \frac{\varepsilon}{2})\frac{2l'}{\varepsilon}$. By (23), we have $\forall l' \geq l_0$,

$$P\left(\exists \underline{\gamma} \in \aleph_{(2l')/\varepsilon}([-2l', 2l']), V(\gamma) \geq \left(\mu - \frac{\varepsilon}{2}\right)\frac{2l'}{\varepsilon}\right) \geq 1 - \frac{\varepsilon^2}{50}.$$

This by translation invariance $\{X_e : e \in E\}$, yields that with probability at least $1 - \frac{\varepsilon}{10}$, for all of the (no more than) $\frac{4}{\varepsilon}$ points

$$\underline{y} \in \Xi_{2l'}\left(\left[-\frac{l'}{\varepsilon} + l', \frac{l'}{\varepsilon} - l'\right] \cap \mathbb{Z}l'\right)$$

with $\underline{y} = (y, 2l')$, there is a path $\underline{\gamma}_y \in \aleph_{((2l')/\varepsilon)+2l'}^y([y - 2l', y + 2l'])$ with:

(i) $|\gamma_y(i) - y| \leq 2l', \forall i \in [0, \frac{2l'}{\varepsilon}]$,
(ii) $V(\gamma_y) \geq (\mu - \frac{\varepsilon}{2})\frac{2l'}{\varepsilon}$.

Now given any $\underline{x} \in \Xi_0([-\frac{l'}{\varepsilon} - l', \frac{l'}{\varepsilon} + l'])$ there exists

$$\underline{y} \in \Xi_{2l'}\left(\left[-\frac{l'}{\varepsilon} + l', \frac{l'}{\varepsilon} - l'\right] \cap \mathbb{Z}l'\right)$$

with $|y - x| \leq 2l'$. For each such $\underline{x}$, pick (arbitrarily) a nonrandom path $\underline{\gamma}_x^1 \in \aleph_{2l'}^x([-\frac{l'}{\varepsilon} - l', \frac{l'}{\varepsilon} + l'])$ with $\gamma_x^1(2l') = y$. For these $\underline{x}$, denote the concatenation of $\underline{\gamma}_x^1$ and $\underline{\gamma}_y$ by $\underline{\gamma}_x$. We then have
$$V(\gamma_x) \geq \left(\mu - \frac{\varepsilon}{2}\right)\frac{2l'}{\varepsilon} + Z,$$
where $Z = \min_x V(\gamma_x^1)$ and each path satisfies $\underline{\gamma}_x \subset [-\frac{l'}{\varepsilon} - l', \frac{l'}{\varepsilon} + l']$. If $Z \geq -l'$, then
$$V(\gamma_x) \geq (\mu - (\mu+1)\varepsilon)\left(\frac{2l'}{\varepsilon} + 2l'\right).$$
But by standard estimates for sums of i.i.d. random variables possessing exponential moments, for some $c > 0$,
$$P(Z \geq -l') \geq 1 - 4\frac{l'}{\varepsilon}e^{-cl'}.$$
Thus,
$$P\left(V(\gamma_x) \geq (\mu - (\mu+1)\varepsilon)\left(\frac{2l'}{\varepsilon} + 2l'\right)\right) \geq 1 - \frac{\varepsilon}{10} - 4\frac{l'}{\varepsilon}e^{-cl'}.$$
Now shift the interval $[-\frac{l'}{\varepsilon} - l', \frac{l'}{\varepsilon} + l']$ to the right by $\frac{l'}{\varepsilon} + l'$ and set $l = (\frac{1}{\varepsilon} + 1)2l'$, relabel $(\mu + 1)\varepsilon$ as $\varepsilon$ and the result holds by increasing $l_0$ if necessary. □

We need a (crude) bound on the lower tail of the distribution of the random variable
$$\min_{\underline{x} \in \Xi_0([0,l])} \max_{\underline{\gamma} \in \aleph_l^{\underline{x}}([0,l))} V(\gamma) = Y_l.$$

The following will suffice

LEMMA 2.2.  *There exists $c > 0$ and $r_0 < \infty$ such that for all $l \geq l_0$*
$$P(Y_l \leq -rl) \leq le^{-crl} \qquad \text{for all } r > r_0.$$

PROOF.  For each $\underline{x} \in \Xi_0([0,l))$ and an arbitrary path $\underline{\gamma}_{\underline{x}} \in \aleph_l^{\underline{x}}([0,l))$, one simply considers the random variable
$$W_{\underline{x}} = \sum_{e \in \gamma_{\underline{x}}} X_e.$$
Recall that $E[X_e] = 0$ and so $E[W_{\underline{x}}] = 0$ as well. By assumption (7), there exists a small positive $c' < c_0$ so that if $e^a = E[e^{-c'X_e}]$, then for $r > r_0$,
$$\begin{aligned}P(W_{\underline{x}} \leq -rl) &= P(-c'W_{\underline{x}} \geq c'rl) \\ &\leq e^{-c'rl}E[e^{-c'W_{\underline{x}}}] \\ &= e^{-c'rl}E[e^{-c'X_e}]^l \\ &\leq e^{-(c'-a/r_0)rl}.\end{aligned}$$

Taking $c'$ small and then $r_0$ large, we find there is a $c > 0$ such that

$$P(W_{\underline{x}} \leq -rl) \leq e^{-crl}, \qquad r > r_0.$$

Since

$$W_{\underline{x}} \leq \max_{\gamma \in \aleph_l^{\underline{x}}([0,l))} V(\gamma),$$

one sees

$$P\left(\max_{\gamma \in \aleph_l^{\underline{x}}([0,l))} V(\gamma) \leq -rl\right) \leq l e^{-crl},$$

from which the result follows. □

In the next proof, we need the following definition.

DEFINITION 2.2. Given $l \in \mathbb{Z}_+$ and $(i,r) \in \{0,1,\ldots,n/l-1\} \times \{0,1,\ldots,n/l-1\}$ say a block $\Xi_{rl}([li, l(i+1)))$ is *good* if for all $\underline{x} \in \Xi_{rl}([li, l(i+1)))$ there exists a path $\underline{\gamma}_{\underline{x}} \in \aleph_{(r+1)l}^{\underline{x}}([li, l(i+1)))$ such that $V(\gamma_{\underline{x}}) \geq (\mu - \varepsilon)l$.

PROOF OF PROPOSITION 2.1. We continue to give proofs in the case $d = 1$. The idea of the proof is to use Lemma 2.1 to show there are on the order of $n$ channels of width $l$ starting at $\Xi_{\delta n}$ and ending at $\Xi_n$. There is a high probability that each channel contains a path with value near $\mu(1-\delta)n$. Then we exploit the independence of the field in these channels.

Given $\varepsilon \in (0, \mu)$, by Lemma 2.1 we can fix $l$ so large that $P(A_l(\varepsilon)) > 1 - c\varepsilon$. Suppose also that $\delta \in (0, 1)$ is given. Without loss of generality, we take $n$ and $(1-\delta)n$ to be multiples of $2l$.

Denote the set of *good* blocks by $\mathcal{G}$ and notice that by our choice of $l$, the random variables $1_{\mathcal{G}}(\Xi_{rl}([li, l(i+1))))$ are i.i.d. Bernoulli random variables with parameter not less than $1 - c\varepsilon$. We partition

$$\Xi_{[\delta n]}((-\delta n/4, \delta n/4)) = \bigcup_{k=1}^{R} \Xi_{[\delta n]}(C_k)$$

into $R = \frac{\delta n}{2l} - 1$ disjoint blocks with $C_k, k = 1, 2, \ldots, R$, of side length $l$ (plus a "remainder" interval) as well as

$$\Xi_{[\delta n, n]}(C_k) = \bigcup_{j=1}^{((1-\delta)n)/l} \Xi_{[\delta n + (j-1)l, \delta n + jl)}(C_k)$$

into $\frac{(1-\delta)n}{l}$ disjoint blocks of side length $l$. Abbreviate the notation by writing

$$R_{k,j} = \Xi_{[\delta n + (j-1)l, \delta n + jl)}(C_k).$$

Fix $k$ and let $Y_{k,j}$ be defined by

(24) $$Y_{k,j} = \min_{\underline{x} \in \Xi_{\delta n + (j-1)l}(C_k)} \max_{\underline{\gamma} \in \aleph^{\underline{x}}_{\delta n + jl}(C_k)} V(\underline{\gamma}).$$

The $Y_{k,j}$ are i.i.d. with lower tail behavior governed by Lemma 2.2. For $c_1$ a small constant, set

(25) $$A(c_1, n, l, k) = \left\{ \exists J \subseteq \left\{ 1, 2, \ldots, \frac{(1-\delta)n}{l} \right\}, |J| \leq \frac{c_1 n}{l}, \right.$$
$$\left. \sum_{j \in J} Y_{k,j} \leq -\frac{\varepsilon}{10}(1-\delta)n \right\}.$$

LEMMA 2.3. *There exists $c'' > 0$ so that for $c_1$ small and $n, l$ sufficiently large, for any $k$ one has*
$$P(A(c_1, n, l, k)) \leq e^{-c'' \varepsilon (1-\delta) n}.$$

PROOF. The statement holds for any $k$ once it holds for one value of $k$ by the translation invariance of the model. Note that by large deviation estimates for binomial random variables,
$$\sum_{j=0}^{c_1 n/l} \binom{(1-\delta)\frac{n}{l}}{j} \leq e^{-(1-\delta) I(c_1/(1-\delta))(n/l)},$$
where $I(\theta) = -\theta \ln \theta - (1-\theta) \ln(1-\theta)$. This bounds the number of subsets $J$ under consideration.

Let $J$ be a subset as described in (25). Using Lemma 2.2 and the constants $c$ and $r_0$ there and Chebyshev bounds, for $c > c' > 0$ not depending on $\varepsilon, l$ or $n$, we have

$$P\left( \sum_{j \in J} Y_{k,j} \leq -\frac{\varepsilon}{10}(1-\delta)n \right)$$
$$\leq e^{-(c'\varepsilon/10)(1-\delta)n} (E[e^{-c'Y}])^{|J|}$$
$$= e^{-(c'\varepsilon/10)(1-\delta)n} (E[e^{-c'Y}; Y \geq -r_0 l] + E[e^{-c'Y}; Y \leq -r_0 l])^{|J|}$$

(26) $$\leq e^{-(c'\varepsilon/10)(1-\delta)n} \left( e^{c'r_0 l} + \frac{l}{c'} \int_{r_0}^{\infty} e^{c'rl} P(Y \leq -rl) \, dr \right)^{(c_1 n)/l}$$

$$\leq e^{-(c'\varepsilon/10)(1-\delta)n} \left( e^{c'r_0 l} + \frac{l}{c'} \int_{r_0}^{\infty} e^{-(c-c')rl} \, dr \right)^{(c_1 n)/l}$$

$$\leq e^{-(c'\varepsilon/10)(1-\delta)n} \left( e^{c'r_0 l} + \frac{l^2}{c'(c-c')} e^{-(c-c')r_0 l} \right)^{(c_1 n)/l}$$

$$\leq e^{-(c\varepsilon/20)(1-\delta)n}, \quad \text{if } c_1 \text{ is small enough.}$$

Thus, $P(A(c_1, n, l, k)) \leq e^{-(I(c_1/1-\delta)(1/l)+c\varepsilon/20)(1-\delta)n}$. Again, by taking $l$ to be large, we obtain the desired bound. $\square$

Now returning to the proof of Proposition 2.1 consider, with fixed $k$,

$$V_k = \sum_{j=1}^{(1-\delta)n/l} 1_{\mathcal{G}}(\Xi_{\delta n+(j-1)l}(C_k)).$$

We have by Lemma 2.3 that $V_k$ is stochastically larger than a binomial random variable with parameters $\frac{(1-\delta)n}{l}$ and $1 - c\varepsilon$. On the event

$$F_k = \left\{V_k \geq \left(\frac{(1-\delta)}{l} - \frac{c_1}{l}\right)n\right\} \cap A(c_1, n, l, k)^c$$

for each point $\underline{x} \in \Xi_{[\delta n]}(C_k)$, there exists a path $\underline{\gamma}_{\underline{x}} \in \aleph_n^{\underline{x}}(C_k)$ which is constructed by concatenating paths with values exceeding $(\mu - \varepsilon/100)l$ through the *good* blocks $R_{k,j}$. When first encountering a *bad* block, select a path from the end point of the path through the *good* block with a value which beats the minimax at (24) in that block and continue in this way through, however, many *bad* blocks it takes until connecting to another *good* block. Upon arriving in a good block, we know we can start from any point and find a path with value exceeding $(\mu - \varepsilon/100)l$ through that block. This gives a connected path. Since we are on the event $A(c_1, n, l, k)^c$, we find that the value of such a concatenated path satisfies

$$\begin{aligned} V(\gamma_{\underline{x}}) &\geq V_k l\left(\mu - \frac{\varepsilon}{100}\right) - (1-\delta)\frac{\varepsilon}{10}n \\ &\geq \left((1-\delta)\left(\mu - \frac{\varepsilon}{100}\right) - c_1\left(\mu - \frac{\varepsilon}{100}\right) - \frac{\varepsilon(1-\delta)}{10}\right)n \\ &\geq (\mu - \varepsilon/5)(1-\delta)n \quad \text{for } c_1 \text{ small enough.} \end{aligned}$$

Let the event

$$D_k = \{\forall \underline{x} \in \Xi_{[\delta n]}(C_k), \exists \underline{\gamma}_{\underline{x}} \in \aleph_n^{\underline{x}}(C_k), V(\gamma_{\underline{x}}) \geq (\mu - \varepsilon/5)(1-\delta)n\}.$$

We have just shown that

$$F_k \subset D_k.$$

Thus, using elementary bounds for the binomial random variables $V_k$ and Lemma 2.3, it follows that

$$P(D_k) \geq 1 - e^{-c_2 \varepsilon n}.$$

Therefore, with $R = \frac{\delta n}{2l}$ and using the independence of the $D_k$, with a new value of $c$ we get

$$P\left(\sum_{k=1}^{R} I_{D_k^c} \leq \frac{R}{10}\right) \geq 1 - K(e^{-c\varepsilon n/2})^{R/10}$$

$$= 1 - e^{-c\varepsilon n^2}.$$

But the discussion above shows that $\{\sum_{k=1}^{R} I_{D_k^c} \leq \frac{R}{10}\}$ is a subset of

$$\{|\{\underline{x} \in \Xi_{[\delta n]}((-\delta n/4, \delta n/4)) : \exists \gamma_{\underline{x}} \in \aleph_n^{\underline{x}},$$

$$V(\gamma_{\underline{x}}) \geq (\mu - \varepsilon/5)(1-\delta)n\}| \geq \tfrac{9}{10}(\delta n/2)\}$$

and Proposition 2.1 is proven. □

We now state the variant of Proposition 2.1 for Model 2. The proof is much the same and so is not given. Here is a rough outline in the case $d = 2$. Outside a set of probability less than $e^{-cn}$, order of $n$ width $l$ channels are proven to exist in which there are paths with values on the order of order $\nu n$. These channels are disjoint and run from $[-\delta n, \delta n]^{d-1} \times \{\delta n\}$ to $[-\delta n, \delta n]^{d-1} \times \{(1-\delta)n\}$. By independence, the probability that none of these channels contains a path with value near $\nu n$ is on the order of $e^{-cn^2}$. Recall the definition of $\Theta_{m,n}(A, B)$ at (9). We now introduce the idea of channels in this context. Set, for $l$ a fixed but large integer, $A_k = [-\delta n + kl, -\delta n + (k+1)l], k = 0, 1, 2, \ldots, \frac{2\delta n}{l}$. Taking $d = 1$ and $B = A_k$, we define the $k$th channel as

$$\Theta_{\delta,n}^k = \{\gamma \in \Theta : \gamma(0) \in A_k \times \{\delta n\}, \gamma(l(\gamma)) \in A_k \times \{(1-\delta)n\},$$

$$\gamma(j) \in A_k \times \mathbf{Z}, j = 0, 1, \ldots, l(\gamma)\}.$$

We then stipulate that a channel $A_k$ is *good* provided

$$\forall (x, y) \in ([-\delta n + kl, -\delta n + (k+1)l] \times \{\delta n\})$$

$$\times ([-\delta n + kl, -\delta n + (k+1)l] \times \{(1-\delta)n\}),$$

$$\exists \gamma \in \Theta_{\delta,n}^k$$

such that $\gamma(0) = x, \gamma(l(\gamma)) = y$ and $V(\gamma) \leq (\nu + \varepsilon/5)(1 - 2\delta)n$

and set

$$\tilde{G}_{\delta n, n} = \left\{|\{k : A_k \text{ is } good\}| \geq \frac{9}{10} \frac{2\delta n}{l}\right\}.$$

The appropriate analog of Proposition 2.1 is the following.

PROPOSITION 2.2. *Consider Model 2. There is an $l_0$ such that for $l \geq l_0$, one has for each sufficiently small $\varepsilon > 0$ and $\delta > 0$ there is a positive $c(\varepsilon, \delta)$ such that*

$$P(\tilde{G}_{\delta n, n}) \geq 1 - e^{-c(\varepsilon, \delta) n^2}.$$

REMARK. For general dimensions $d \geq 2$ with events analogously defined the bound will be $1 - e^{-c(\varepsilon, \delta) n^d}$.

**3. Near the starting point: Gaussian case.** In this section, we prove Theorem 1.2. To simplify the exposition, we consider $n$ of the form $2^N$. The conclusions we arrive at will easily be seen to hold for arbitrary large $n$. We first consider a lower bound for

$$P(Z_{2^N} \leq (\mu - \varepsilon) 2^N).$$

Then we show it is of the correct logarithmic order by considering the upper bound. Recall the notation $\underline{x} = (x, m) \in \Xi_m$ where $x \in \mathbb{Z}$ and $m \in \mathbb{Z}^+$. For $k = 0, 1, \ldots$, define

$$T_k = \{e \in E : e = ((x, m), (x \pm 1, m + 1)), |x| \leq m, m \in [2^k, 2^{k+1})\}.$$

Notice that $|T_k| \leq 2^{2k+1}$. In higher dimensions, we have $|T_k| \leq 2d(22^{k+1} + 1)^d 2^k$. [In order to accomodate the origin $(0,0)$, we include the edges from this point in $T_0$.] In the following, things are described for $T_k$ with $k > 0$; we rely on the reader to make the necessary adjustments to include $T_0$ in the arguments. Fix $M > 0$ and set

$$(27) \qquad A_k^N = \left\{ \forall e \in T_k, X_e \leq -\frac{M 2^{N-k}}{N} \right\}.$$

This event would cause a lower deviation in the value of $Z_n$. Using the fact that the random variables $\{X_e : e \in T_k\}$ are i.i.d. $\mathcal{N}(0,1)$, the following results are easily seen to hold.

LEMMA 3.1. *For $A_k^N$ as above with $\{X_e : e \in T_k\}$ i.i.d. $\mathcal{N}(0,1)$, there exists a positive constant $c$ so that for all $N, M$ and $k \leq \frac{N}{2}$,*

$$P(A_k^N) \geq e^{-c M^2 2^{2N}/N^2}.$$

PROOF. This is a simple combination of the fact that $\{X_e : e \in T_k\}$ are i.i.d. $\mathcal{N}(0,1)$ and that there are no more than $2^{2(k+1)}$ points in $T_k$ with $k \leq \frac{N}{2}$. □

Since the events $\{A_k^N : 0 \leq k \leq \frac{N}{2}\}$ are independent, it immediately follows that

COROLLARY 3.1. *For $A_k^N$ as above, there exists a positive constant $c$ so that*
$$P\left(\bigcap_{k=0}^{N/2} A_k^N\right) \geq e^{-cM^2 2^{2N}/N}.$$

Finally, we have:

LEMMA 3.2. *If*
$$B_{N/2}^N = \{\forall \underline{\gamma} \in \aleph_{2^{N/2+1},2^N}^{[-2^{N/2+1},2^{N/2+1}]}, V(\gamma) \leq (\mu+\varepsilon)2^N\},$$
*then*
$$\lim_{N\to\infty} P(B_{N/2}^N) = 1 \quad a.s.$$

PROOF. By (15), it follows that for any $\underline{x} \in \Xi_{2^{N/2+1}}([-2^{N/2+1},2^{N/2+1}])$, for some $c > 0$,
$$P\left(\sup_{\underline{\gamma} \in \aleph_{2^N}^{\underline{x}}} V(\gamma) > (\mu+\varepsilon)2^N\right) \leq e^{-c2^N}.$$
Since there are at most $2^{N/2+2}+1$ points in $\Xi_{2^{N/2+1}}([-2^{N/2+1},2^{N/2+1}])$, the result follows. □

Putting these together, we can obtain the following proposition.

PROPOSITION 3.1. *Given $\varepsilon > 0$, there is a $c(\varepsilon) > 0$ so that for all $N$ large,*
$$P(Z_{2^N} \leq (\mu-\varepsilon)2^N) \geq 1/2 e^{-c(\varepsilon)2^{2N}/N}.$$

PROOF. By the independence of the field $\{X_e : e \in E\}$, the events $\bigcap_{k=0}^{N/2} A_k^N$ and $B_{N/2}^N$ are independent. Thus, by Lemmas 3.1 and 3.2, for $N$ large
$$P\left(\bigcap_{k=0}^{N/2} A_k^N \cap B_{N/2}^N\right) \geq \tfrac{1}{2} e^{-cM^2 2^{2N}/N}.$$
Now on the event $\bigcap_{n=0}^{N/2} A_n^N \cap B_{N/2}^N$, noting that any path $\underline{\gamma} \in \aleph_{2^N}$ will take $2^k$ steps in $T_k$, if $M \geq 4\varepsilon$
$$V(\gamma) \leq \sum_{k=0}^{N/2} -\frac{M 2^{N-k}}{N} 2^k + (\mu+\varepsilon)2^N$$
(28)
$$< -\frac{M 2^N}{2} + (\mu+\varepsilon)2^N$$
$$\leq (\mu-\varepsilon)2^N.$$

LARGE DEVIATION REGIMES 17

Thus, we have the result with $c(\varepsilon) = 16c\varepsilon^2$. □

Immediately, we get the first inequality in Theorem 1.2, namely, for all $\varepsilon > 0$,

$$\lim_{n\to\infty} \frac{\log n}{n^2} \log P(Z_n \leq (\mu - \varepsilon)n) > -\infty.$$

It remains to show that $\frac{n^2}{\log n}$ does indeed give the correct (logarithmic) rate by showing

$$\overline{\lim_{n\to\infty}} \frac{\log n}{n^2} \log P(Z_n \leq (\mu - \varepsilon)n) < 0.$$

When analyzing $\underline{\lim}_{n\to\infty} \frac{\log n}{n^2} \log P(Z_n \leq (\mu - \varepsilon)n) < 0$, we could simply consider the event where all values of $\{X_e : e \in T_k\}, 0 \leq k \leq \frac{N}{2}$ were large negative values. However, for the case of $\overline{\lim}_{n\to\infty} \frac{\log n}{n^2} \log P(Z_n \leq (\mu - \varepsilon)n) < 0$, we will need to consider events on which the sums of the negative parts of the $X_e$ for $e \in T_k$ take on large values and also incorporate the information from Proposition 2.1. Since the difficulties arise from negative values of $X_e$ we are naturally led to consider, for $0 \leq k \leq N + \log \delta$, the (independent) random variables

$$(29) \qquad V_k^- = \sum_{e \in T_k} X_e^-,$$

where $X_e = X_e^+ - X_e^-$. These will control the negative values of the field near the starting point. We first address bounds on $V_k^-$.

LEMMA 3.3. *Given $\delta > 0$ a negative power of 2, integers $k$ and $N$ with $0 \leq k \leq N + \log \delta$ and $c > 2^{k-N+3}$ we have with $c_0 = 1 - \frac{\log 2}{2}$,*

$$P(2^{-k} V_k^- \geq c2^N) \leq \exp\left\{\frac{-c_0 c^2 2^{2N}}{16}\right\}.$$

PROOF. For the first inequality, we observe that $X_e^- \leq |X_e|$. Thus, for any $\lambda > 0$,

$$P(2^{-k} V_k^- \geq c2^N) \leq e^{-\lambda c 2^N} E[e^{\lambda 2^{-k} V_k^-}]$$

$$= e^{-\lambda c 2^N} E[e^{\lambda 2^{-k} X_e^-}]^{|T_k|}$$

$$\leq e^{-\lambda c 2^N} E[e^{\lambda 2^{-k} |X_e|}]^{|T_k|}$$

$$(30) \qquad = e^{-\lambda c 2^N} \left(e^{(\lambda 2^{-k})^2/2} \sqrt{\frac{2}{\pi}} \int_0^\infty e^{-(x - \lambda 2^{-k})^2/2} \, dx\right)^{|T_k|}$$

$$\leq e^{-\lambda c 2^N} \left( e^{(\lambda 2^{-k})^2/2} \sqrt{\frac{2}{\pi}} \int_{-\lambda 2^{-k}}^{\infty} e^{-x^2/2} \, dx \right)^{|T_k|}$$

$$\leq e^{-\lambda c 2^N} (2 e^{(\lambda 2^{-k})^2/2})^{|T_k|}$$

$$\leq e^{-\lambda c 2^N + 2^{2(k+1)} \log 2 + 2\lambda^2}.$$

Since this holds for all $\lambda > 0$, we are free to select $\lambda = c 2^{N-2}$. Recalling that $c > 2^{k-N+3}$ (30) becomes

$$P(2^{-k} V_k^- \geq c 2^N) \leq e^{-(c^2 2^{2N})/4 + 2^{2(k+1)} \log 2 + (c^2 2^{2N})/8}$$

(31)
$$= e^{2^{2N} - c^2/8 + 2^{2(k-N+1)} \log 2}$$

$$\leq e^{-(c_0 c^2/16) 2^{2N}}. \qquad \square$$

LEMMA 3.4. *There is a $c > 0$ such that given $\varepsilon > 0$ and $\delta > 0$ for which $2^7 \delta < \varepsilon$, we have for all $N$ sufficiently large,*

$$P\left( \sum_{k=0}^{N + \log \delta} 2^{-k} V_k^- \geq \varepsilon 2^N \right) \leq e^{-(c \varepsilon^2 2^{2N})/N}.$$

PROOF. The strategy of the proof is to handle the disorder which can occur subject to the condition $\sum_{k=0}^{N+\log \delta} 2^{-k} V_k^- \geq \varepsilon 2^N$. Denote $\mathcal{Z} = (\mathbb{Z}_+ \cup \{\infty\})^{N + \log \delta + 1}$, where again (and throughout the rest of the paper), we assume without losing generality that $\log \delta$ is an integer. For $\mathbf{v} \in \mathcal{Z}$, with $\mathbf{v} = (v_0, v_1, \ldots, v_{N + \log \delta})$, set

$$\mathcal{I}(\mathbf{v}) = \{k \in \{0, 1, 2, \ldots, N + \log \delta\} : 2^{-v_k} > 2^{k-N+3}\}.$$

Now put

$$\mathcal{V} = \left\{ \mathbf{v} : \sum_{k \in \mathcal{I}(\mathbf{v})} 2^{-v_k} \geq \frac{3\varepsilon}{8} \right\}$$

and for $\mathbf{v} \in \mathcal{V}$, define the event

$$B(\mathbf{v}) = \bigcap_{k \in \mathcal{I}(\mathbf{v})} \{2^{-k} V_k^- \geq 2^{N - v_k}\}.$$

Our goal is to estimate the probability of $A \equiv \{\sum_{k=0}^{N + \log \delta} 2^{-k} V_k^- \geq \varepsilon 2^N\}$ by showing

$$A \subset \bigcup_{\mathbf{v} \in \mathcal{B}} B(\mathbf{v})$$

for $\mathcal{B}$ a suitable subset of $\mathcal{V}$. Now write

$$I(\omega) = \{k \in \{0, 1, 2, \ldots, N + \log \delta\} : 2^{-k} V_k^- \geq 2^{k+3}\}.$$

Our assumption $2^7\delta < \varepsilon$ implies that

(32)
$$\sum_{k \notin I(\omega), k=0}^{N+\log \delta} 2^{-k} V_k^- \leq \sum_{0}^{N+\log \delta} 2^{k+3}$$
$$\leq 2^{N+4+\log \delta}$$
$$= \delta 2^{N+4}$$
$$\leq \frac{\varepsilon}{4} 2^N,$$

by our choice of $\delta$. Thus, if $\omega \in A$,

$$\sum_{k \in I(\omega)} 2^{-k} V_k^- \geq \frac{3\varepsilon}{4} 2^N.$$

Given $\omega \in A$, we now produce a $\mathbf{v}$ for which $\omega \in B(\mathbf{v})$. For $k \in I(\omega)$, and $2^{-k} V_k^- < 2^N$, select $v_k \in \mathbb{Z}_+$ satisfying

$$2^{N-v_k+1} > 2^{-k} V_k^- \geq 2^{N-v_k}.$$

When $k \in I(\omega)$ and $2^{-k} V_k^- \geq 2^N$, take $v_k = 0$. For $k \notin I(\omega)$, select $v_k = \infty$. This gives a $\mathbf{v} \in \mathcal{V}$ for which

(33)
$$\sum_{k \in I(\omega)} 2^{-v_k} \geq \frac{3\varepsilon}{8},$$

since if $2^{-k'} V_{k'}^- \geq 2^N$ for some $k' \in I(\omega)$, then

$$\sum_{k \in I(\omega)} 2^{-v_k} \geq 2^{-v'_k} = 1,$$

while if $2^{-k} V_k^- < 2^N$ for all $k \in I(\omega)$, then

$$\sum_{k \in I(\omega)} 2^{-v_k} \geq 2^{-N-1} \sum_{k \in I(\omega)} 2^{-k} V_k^- \geq \frac{3\varepsilon}{8}.$$

Thus, $\omega \in B(\mathbf{v})$. Denote the set of vectors $\mathbf{v} \in \mathcal{V}$ such that either $2^{-v_k} \geq 2^{k-N+3}$ or $v_k = \infty$ by $\mathcal{B}$. We note that $|\mathcal{B}| \leq N^N$. Then by (32),

(34)
$$A \subset \bigcup_{\mathbf{v} \in \mathcal{B}} B(\mathbf{v}).$$

Moreover, by Lemma 3.3,

$$P(B(\mathbf{v})) \leq \exp\left\{-\frac{c_0 2^{2N}}{16} \sum_{k \in I(\mathbf{v})} 2^{-2v_k}\right\}$$

$$\leq \exp\left\{-\frac{c_0 2^{2N}}{16N}\left(\sum_{k\in I(\mathbf{v})} 2^{-v_k}\right)^2\right\} \quad \text{(by Cauchy–Schwarz)}$$

(35)
$$\leq \exp\left\{-\frac{c_0 2^{2N}}{16N}\left(\frac{3\varepsilon}{8}\right)^2\right\}$$

$$= \exp\left\{-\frac{c\varepsilon^2 2^{2N}}{N}\right\},$$

for $c > 0$ not depending on $n$. To end the proof, by (34) and (35), we have

$$P(A) \leq N^N \exp\left\{-\frac{c\varepsilon^2 2^{2N}}{N}\right\},$$

which gives the claim on adjusting the value of $c$. □

Thus, in order to prove Theorem 1.2, it will suffice to show that given $\varepsilon > 0$ and a constant $\eta$ to be named later, there is a $c_1 = c(\eta, \varepsilon)$ such that

$$(36) \quad P\left(\{Z_{2^N} \leq (\mu - \varepsilon)2^N\} \cap \left\{\sum_{k=0}^{N+\log\delta} 2^{-k}V_k^- \leq \eta\varepsilon 2^N\right\}\right) \leq e^{-c_1 2^{2N}/N}.$$

In fact, given Lemma 3.4, this will follow from the stronger inequality

$$(37) \quad P\left(\{Z_{2^N} \leq (\mu - \varepsilon)2^N\} \cap \left\{\sum_{k=0}^{N+\log\delta} 2^{-k}V_k^- \leq \eta\varepsilon 2^N\right\}\right) \leq e^{-c 2^{2N}}$$

for some strictly positive $c$. We pick a random path $\underline{\gamma}^k \in \aleph_{2^k, 2^{k+1}}$ with

$$\underline{\gamma}^k(0) \in \Xi_{2^k}\left(\left[\frac{-2^k}{4}, \frac{2^k}{4}\right]\right),$$

$$\underline{\gamma}^k(2^{k+1} - 2^k) = \gamma^k(2^k) \in \Xi_{2^{k+1}}\left(\left[\frac{-2^{k+1}}{4}, \frac{2^{k+1}}{4}\right]\right)$$

as follows:

1. Independently of the field $\{X_e : e \in E\}$, pick the initial site $\underline{\gamma}^k(0)$ uniformly in $\Xi_{2^k}([\frac{-2^k}{4}, \frac{2^k}{4}])$.

2. Independently of the field $\{X_e : e \in E\}$, and the initial point, pick the terminal site $\underline{\gamma}^k(2^k)$ uniformly in $\Xi_{2^{k+1}}([\frac{-2^{k+1}}{4}, \frac{2^{k+1}}{4}])$.

3. Given $\underline{\gamma}^k(i)$, $\underline{\gamma}^k(i+1)$ is deterministically fixed to be the nearest neighbor closest to $\underline{\gamma}^k(2^k)$. If there are two equal closest next moves, it moves to the left. We shall denote probabilities and expectations with respect to this random selection procedure by $\tilde{P}$ and $\tilde{E}$.

The proof of the following lemma is left to the reader.

LEMMA 3.5. *There is a positive c so that for all k and for every edge $e \in T_k$,*

$$\tilde{P}(e \in \underline{\gamma}^k) \leq c2^{-k}.$$

COROLLARY 3.2. *There is a positive c so that for k a strictly positive integer and $\underline{\gamma}^k$, as above, we have*

$$\tilde{P}\left(\sum_{e \in \underline{\gamma}^k} X_e^- \geq 100c2^{-k}V_k^-\right) \leq \tfrac{1}{100}.$$

PROOF. We have

$$\tilde{E}\left[\sum_{e \in \underline{\gamma}^k} X_e^-\right] = \sum_{e \in T_k} \tilde{P}(e \in \underline{\gamma}^k)X_e^-$$

$$\leq c2^{-k}\sum_{e \in T_k} X_e^-, \qquad \text{by Lemma 3.5}$$

$$= c2^{-k}V_k^-,$$

so the result is simply the Markov inequality. □

The following is also a consequence of Markov's inequality and is left for the reader to prove.

COROLLARY 3.3. *For $\underline{x} \in \Xi_{2^k}([\tfrac{-2^k}{4}, \tfrac{2^k}{4}])$, let $\Gamma(\underline{x})$ represent the points $\underline{y}$ in $\Xi_{2^{k+1}}([\tfrac{-2^{k+1}}{4}, \tfrac{2^{k+1}}{4}])$ such that there exists a path from $\underline{x}$ to $\underline{y}$, $\underline{\gamma}_{\underline{x},\underline{y}}$, for which $\sum_{e \in \underline{\gamma}_{\underline{x},\underline{y}}} X_e^- \leq 100c2^{-k}V_k^-$. For $\tfrac{9}{10}$ of the points $\underline{x} \in \Xi_{2^k}([\tfrac{-2^k}{4}, \tfrac{2^k}{4}])$, the cardinality of $\Gamma(\underline{x})$ is at least $\tfrac{9}{10}$ that of $\Xi_{2^{k+1}}([\tfrac{-2^{k+1}}{4}, \tfrac{2^{k+1}}{4}])$. Equally for $\tfrac{9}{10}$ of the points $\underline{y} \in \Xi_{2^{k+1}}([\tfrac{-2^{k+1}}{4}, \tfrac{2^{k+1}}{4}])$, the cardinality of the set $\{\underline{x} \in \Xi_{2^k}([\tfrac{-2^k}{4}, \tfrac{2^k}{4}]) : \underline{y} \in \Gamma(\underline{x})\}$ is at least $\tfrac{9}{10}$ that of $\Xi_{2^k}([\tfrac{-2^k}{4}, \tfrac{2^{k+1}}{4}])$.*

Assume without loss of generality that $100c > 1$ and select $\eta = \tfrac{1}{200c}$ and prove by induction.

COROLLARY 3.4. *For all $k = 0, 1, 2, \ldots, N + \log \delta$ there is a path, $\underline{\gamma}_{\underline{y}}$, from $\underline{0}$ to $\tfrac{9}{10}$ of the points $\underline{y} \in \Xi_{2^{k+1}}([-\tfrac{2^{k+1}}{4}, \tfrac{2^{k+1}}{4}])$, such that*

$$\sum_{e \in \gamma_{\underline{y}}} X_e^- \leq 100c \sum_{i=0}^{k} 2^{-i}V_i^-.$$

*In particular, on the event* $\{\sum_{k=0}^{N+\log \delta} 2^{-k} V_k^- \leq \eta \varepsilon 2^N\}$, *for* $\frac{9}{10}$ *of the points,* $\underline{y} \in \Xi_{\delta 2^N}([-\frac{\delta 2^N}{4}, \frac{\delta 2^N}{4}])$, *there is a path* $\underline{\gamma}_{\underline{y}} \in \aleph_{\delta 2^N}$ *such that*

$$V(\gamma_{\underline{y}}) \geq -\frac{\varepsilon}{2} 2^N. \tag{38}$$

PROOF. By induction, it is plainly true for $k = 0$ as we simply consider a single random variable. Now given that it is true for $k$ we have by Corollary 3.3 that for $\frac{9}{10}$ of $\underline{y} \in \Xi_{2^{k+2}}([-\frac{2^{k+2}}{4}, \frac{2^{k+2}}{4}])$, there is a path from $\underline{x}$ to $\underline{y}$, $\underline{\gamma}_{\underline{x},\underline{y}}$, so that for $\frac{9}{10}$ of the $\underline{x} \in \Xi_{2^{k+1}}([-\frac{2^{k+1}}{4}, \frac{2^{k+1}}{4}])$, we have

$$\sum_{e \in \underline{\gamma}_{\underline{x},\underline{y}}} X_e^- \leq 100c 2^{-(k+1)} V_{k+1}^-.$$

But as $\frac{9}{10} + \frac{9}{10} > 1$, it must be the case that for such a $\underline{y}$ there exists such an $\underline{x}$ also with the property that there exists a $\underline{\gamma}_{\underline{x}}$ from $\underline{0}$ to $\underline{x}$ such that

$$\sum_{e \in \underline{\gamma}_{\underline{x}}} X_e^- \leq 100c \sum_{i=0}^{k} 2^{-i} V_i^-.$$

The path $\underline{\gamma}_{\underline{y}}$ obtained by the concatenation of $\underline{\gamma}_{\underline{x}}$ and $\underline{\gamma}_{\underline{x},\underline{y}}$ gives a path from $\underline{0}$ to $\underline{y}$ with the desired properties. □

PROOF OF THEOREM 1.2. Define

$$F_{N,\delta} = \left\{ \text{for } \frac{9}{10} \frac{2^N \delta}{4} \text{ sites } \underline{y} \in \Xi_{\delta 2^N}\left([-\frac{2^N \delta}{4}, \frac{2^N \delta}{4}]\right), \right.$$

$$\left. \forall \gamma \in \aleph_{2^N}^{\underline{y}}, \sum_{e \in \gamma} X_e \leq \left(\mu - \frac{\varepsilon}{2}\right) 2^N \right\}.$$

By the construction of the path in Corollary 3.4,

$$\{Z_{2^N} \leq (\mu - \varepsilon) 2^N\} \cap \left\{ \sum_{k=0}^{N+\log \delta} 2^{-k} V_k^- \leq \eta \varepsilon 2^N \right\} \subset F_{N,\delta}.$$

But by Proposition 2.1, for $\delta$ sufficiently small and $\varepsilon \ll 1$,

$$P(F_{N,\delta}) \leq e^{-c' 2^{2N}}.$$

Thus, by Lemma 3.4,

$$P(Z_{2^N} \leq (\mu - \varepsilon) 2^N) \leq P\left(\bigcup_{\mathbf{v}} B(\mathbf{v})\right) + P(F_{N,\delta})$$

$$\leq e^{-c 2^{2N}/N} + e^{-c' 2^{2N}}. \tag{39}$$

□

**4. Near the starting point: sub-Gaussian case.** We continue to consider the oriented percolation model on $\Xi$ where the directed edges go from a site $(x,n)$ to the sites $(x \pm 1, n+1)$, in $d = 1$. The development in higher dimensions is analogous. Recall the $\{X_e : e \in E\}$ are i.i.d. of mean 0 and so interpreted as passage times, are not necessarily positive. Also, recall that $X_e^-$ denotes the negative part of $X_e$. We assume that there exists a positive constant $M_0$ and a increasing function $f$ so that

(40) $\quad P(X_e < -x) = P(X_e^- > x) = e^{-x^2 f(x)} \qquad \text{for } x \geq M_0, \ f(M_0) > 0.$

We now prove Theorem 1.1, that is, for all $\varepsilon > 0$, there exists a positive constant $C = C(\varepsilon)$ such that

$$\overline{\lim_{n \to \infty}} \frac{1}{n^2} \log P\left(\sup_{\gamma \in \aleph_n} V(\gamma) \leq (\mu - \varepsilon)n\right) \leq -C$$

if and only if

(41) $$\sum_{k=1}^{\infty} \frac{1}{f(2^k)} < \infty.$$

We prove the result for $n$ of the form $2^N$, but it will be clear from the proof that the approach extends to general $n$. In evaluating the justness of the assumptions on the distribution of the $X_e$, it is worth remarking that the work of the preceding section can be adapted to show that if the distribution of the $X_e$ satisfies

$$\liminf_{x \to \infty} -\frac{\log P(X_e < x)}{x^2} < \infty,$$

then the limsup of Theorem 1.1 will equal 0.

Recall that for $k \geq 1, \ldots, T_k$ denotes the set of edges $e$ from sites $(x, m)$ with $|x| \leq 2^{k+1}$ and $2^k \leq m < 2^{k+1}$. We abbreviate the notation by writing

$$R_{2^k} = \Xi_{2^k}([-2^k, 2^k]).$$

Let $\delta > 0$ and assume that $\log \delta$ is an integer and that $N + \log \delta \gg 0$. Recalling the definition of $V_k^-$ from (29), define the events, for $k = 1, \ldots, N + \log \delta$,

(42) $\qquad B_k(v) = \{2^{-k} V_k^- \in [2^{N-v}, 2^{N-v+1})\}$

and, for an $N$-tuple $\mathbf{v} = (v_1, \ldots, v_N)$ of nonnegative integers or infinity, the event

(43) $$A(\mathbf{v}) = \bigcap_{k=1}^{N+\log \delta} B_k(v_k).$$

More generally, given $I \subset \{1, 2, \ldots, N + \log \delta\}$ and $\mathbf{v} = (v_1, \ldots, v_{N+\log \delta})$, set

(44) $$A^I(\mathbf{v}) = \bigcap_{k \in I} B_k(v_k).$$

Let $M$ be large enough so that $M > M_0$, with $M_0$ appearing in (40). For each edge $e$, let
$$X_e^{-,M} = X_e^- 1_{\{X_e^- > M\}}.$$
As is the case for $X_e^-$, we also have $P(X_e^{-,M} > x) \leq e^{-x^2 f(x)}$ for $x \geq M_0$. In order to prove Theorem 1.1, we use

PROPOSITION 4.1. *Under assumptions (40) and (41), for $\varepsilon > 0$, there is an $\varepsilon_1 = \varepsilon_1(\varepsilon) > 0$ and a $\delta$ (of the form $2^{-l}$) for integer $l$, so that*
$$\sum_{k=-\log \delta}^{\infty} 2^{-k/2} \leq \frac{\varepsilon_1}{8}$$
*and so that, for $I_{N+\log \delta} = \{\mathbf{v} = (v_1, \ldots, v_{N+\log \delta}) \in \mathbb{Z}_+^{N+\log \delta} : \sum_{k=1}^{N+\log \delta} 2^{-v_k} \geq \varepsilon_1\}$:*

1. $P(\bigcup_{\mathbf{v} \in I_{N+\log \delta}} A(\mathbf{v})) \leq e^{-C 2^{2N}}$ *for some $C = C(\varepsilon)$ not depending on $N$;*
2. *for $\mathbf{v} \notin I_{N+\log \delta}$, on the event $A(\mathbf{v})$ for $\frac{9}{10}$ of the points $x$ in the middle quarter of $R_{2^{N+\log \delta}}$, there is at least one path $\gamma_0^x$ from $(\underline{0},0)$ to $x$ so that $V(\gamma_0^x) \geq -\frac{\varepsilon}{2} 2^N$.*

We record the next result which is an immediate consequence of the hypothesis (40).

LEMMA 4.1. *For $M > M_0$, let $Y$ be a random variable with a $\mathcal{N}(0, \frac{1}{f(M)})$ distribution and set $Y^M = Y 1_{\{Y > M\}}$. Then $X_e^{-,M}$ is stochastically less than $Y^M$, that is, to say*
$$P(X_e^{-,M} > x) \leq P(Y^M > x) \qquad \forall x \in \mathbb{R}.$$

Building on Lemma 4.1, we have the following lemma.

LEMMA 4.2. *There exists universal positive $C$ so that for any integer $v$ with $2^{N-k-v-2} = M \geq M_0$,*

(45)
$$P(B_k(v)) \leq e^{-C 2^{2(N-v)} f(M)}.$$

PROOF. Since $M = 2^{N-v-k-2}$ and $|T_k| \leq 2^{2k+1}$,
$$\sum_{e \in T_k} X_e^- - \sum_{e \in T_k} X_e^{-,M} = \sum_{e \in T_k} X_e^- 1_{\{X_e^- \leq M\}}$$
(46)
$$\leq 2^{2k+1} M$$
$$= 2^{N-v+k-1}.$$

Therefore, as $\sum_{e \in T_k} X_e^- \geq 2^{N-v+k}$ on $B_k(v)$, we have

$$B_k(v) \subset \left\{ \sum_{e \in T_k} X_e^{-,M} \geq 2^{N-v+k-1} \right\}.$$

On the other hand, for $e \in T_k$, let $\{Y_e : e \in E\}$ be i.i.d. $\mathcal{N}(0, 1/f(M))$ random variables and $Y_e^M = Y_e 1_{\{Y_e > M\}}$. We have for $\theta = M f(M)$,

$$E[e^{\theta Y_e^M}] = P(Y_e \leq M) + \sqrt{\frac{f(M)}{2\pi}} \int_M^\infty \exp\left\{ \theta y - \frac{y^2 f(M)}{2} \right\} dy$$

$$= P(Y_e \leq M) + \sqrt{\frac{f(M)}{2\pi}} \int_M^\infty \exp\left\{ \frac{\theta^2}{2f(M)} - \left[ y - \frac{\theta}{f(M)} \right]^2 \frac{f(M)}{2} \right\} dy$$

$$= P(Y_e \leq M) + \frac{1}{2} \exp\left\{ \frac{\theta^2}{2f(M)} \right\}.$$

As $P(Y_e \leq M) < 1 < \frac{1}{2}e$, we have

$$E[e^{\theta Y_e^M}] \leq \exp\left\{ \frac{\theta^2}{2f(M)} \right\}$$

for $\theta^2 = M^2 f^2(M) > 2 f(M)$, which holds for $M$ large. Thus, by Lemma 4.1 and the Chebyshev inequality,

$$P(B_k(v)) \leq P\left( \sum_{e \in T_k} X_e^{-,M} \geq 2^{N-v+k-1} \right)$$

$$\leq \exp\{-\theta 2^{N-v+k-1}\} (E[e^{\theta X_e^{-,M}}])^{|T_k|}$$

$$\leq \exp\{-\theta 2^{N-v+k-1}\} (E[e^{\theta Y_e^M}])^{|T_k|}$$

$$\leq \exp\{-\theta 2^{N-v+k-1}\} (E[e^{\theta Y_e^M}])^{2^{2k+1}}$$

$$\leq \exp\left\{ \frac{\theta^2}{2f(M)} 2^{2k+1} - \theta 2^{N-v+k-1} \right\}$$

$$= \exp\{f(M) 2^{2(N-v)-4} - f(M) 2^{2(N-v)-3}\}$$

$$\leq \exp\left\{ -\frac{1}{16} f(M) 2^{2(N-v)} \right\}$$

and the lemma is proved with $C = \frac{1}{16}$. $\square$

We now prove claim 1 of Proposition 4.1.

PROOF OF PROPOSITION 4.1. We must show

$$P\left(\bigcup_{\mathbf{v}\in I_{N+\log\delta}} \bigcap_{k=1}^{N+\log\delta} B_k(v_k)\right) \leq e^{-C 2^{2N}}.$$

For $\mathbf{v} = (v_1, \ldots, v_{N+\log\delta}) \in I_{N+\log\delta}$, we say

$$v_k \text{ is } bad, \quad \text{if } v_k > \frac{N-k}{2} - 2$$

or equivalently we say that

$$v_k \text{ is } good, \quad \text{if } v_k \leq \frac{N-k}{2} - 2.$$

Let $G = \{k_1, \ldots, k_j\}, j \leq N + \log\delta$, be the indices corresponding to *good* $v_k$'s and $B = \{k_{j+1}, \ldots, k_{N+\log\delta}\}$ be the indices corresponding to *bad* ones. Define now the event

$$A^G(\mathbf{v}) = \bigcap_{i=1}^{j} B_{k_i}(v_{k_i})$$

and note that there are at most $2^{N+\log\delta}$ possible choices of $G$ and that for a given choice of $k_1, \ldots, k_j$, there are at most $\frac{N}{2}$ possible choices for each of the *good* $v_{k_i}$'s. Hence, if it is proven that (for some $C$ not depending on the particular $v$)

$$P(A^G(\mathbf{v})) \leq e^{-C 2^{2N}},$$

then it follows that for $N$ large enough,

$$P\left(\bigcup_{\mathbf{v}\in I_{N+\log\delta}} A(\mathbf{v})\right) \leq 2^N \left(\frac{N}{2}\right)^N e^{-C 2^{2N}}$$

(47)
$$\leq e^{-C 2^{2N}} \quad \text{with a new choice of } C.$$

Notice that Lemma 4.2 applies for $v = v_k$ *good* provided $\frac{1}{\delta} > 4M_0$. Thus, using the assumption that $f$ is nondecreasing, we have

$$P(A^G(\mathbf{v})) = \prod_{i=1}^{j} P(B_{k_i}(v_{k_i}))$$

(48)
$$\leq \exp\left\{-C 2^{2N} \sum_{i=1}^{j} 2^{-2v_{k_i}} f(2^{N-v_{k_i}-k_i-2})\right\}.$$

Notice that $2^{-v_{k_i}} < 4 \cdot 2^{-(N-k_i)/2}$ for $i \in B = \{j+1, \ldots, N + \log \delta\}$, so their contribution is

$$4 \sum_{i=j+1}^{N+\log \delta} 2^{-v_{k_i}} < 4 \sum_{i=j+1}^{N+\log \delta} 2^{-(N-k_i)/2}$$

$$\leq 4 \sum_{k=1}^{N+\log \delta} 2^{-(N-k)/2}$$

$$\leq 4 \sum_{k=\log(1/\delta)}^{N-1} 2^{-k/2}$$

$$\leq 4 \sum_{k=\log(1/\delta)}^{\infty} 2^{-k/2}$$

$$\leq \frac{\varepsilon_1}{2}.$$

Then as $\sum_{k=1}^{N+\log \delta} 2^{-v_k} \geq \varepsilon_1$,

$$\tag{49} \sum_{i=1}^{j} 2^{-v_{k_i}} \geq \frac{\varepsilon_1}{2}.$$

We use this to show that under the condition $\sum_{k \in G} 2^{-v_k} \geq \frac{\varepsilon_1}{2}$ there is a positive $C$ such that

$$\sum_{k \in G} 2^{-2v_k} f(2^{N-v_k-k-2}) \geq \sum_{k \in G} 2^{-2v_k} f(2^{(N-k)/2})$$

$$\geq C.$$

From the Cauchy–Schwarz inequality together with hypothesis (49),

$$\tag{50} \frac{\varepsilon_1^2}{(2 \sum_{k \in G} 1/f(2^{(N-k)/2}))^2} \leq \frac{\sum_{k \in G} 2^{-2v_k} f(2^{(N-k)/2})}{\sum_{k \in G} 1/f(2^{(N-k)/2})},$$

that is,

$$\tag{51} \frac{\varepsilon_1^2}{4 \sum_{k \in G} 1/f(2^{(N-k)/2})} \leq \sum_{k \in G} 2^{-2v_k} f(2^{(N-k)/2}).$$

Since $f$ is nondecreasing,

$$0 < \sum_{k \in G} 1/f(2^{(N-k)/2})$$

$$\leq \sum_{k=1}^{N+\log \delta} 1/f(2^{(N-k)/2})$$

$$= \sum_{k=-\log \delta}^{N-1} 1/f(2^{k/2})$$

(52)
$$\leq \sum_{k=2}^{\infty} 1/f(2^{k/2})$$

$$\leq 2 \sum_{k=1}^{\infty} 1/f(2^k)$$

$$< \infty.$$

So, by (51) and (52), we have

$$\sum_{k=1}^{N} 2^{-2v_k} f(2^{(N-k)/2}) \geq \frac{\varepsilon_1^2}{4 \sum_{k \in G} 1/f(2^{(N-k)/2})}$$
$$> 0.$$

This proves part one of Proposition 4.1. □

REMARK. For the higher dimensional cases, the major difference in proof is that we must use Holder's inequality rather than the Cauchy–Schwarz.

We now turn to the proof of claim 2 of Proposition 4.1. This involves a continuation of the ideas of Lemma 3.3.

PROOF. Suppose that $A(\mathbf{v})$ occurs for some $\mathbf{v} = (v_1, \ldots, v_{N+\log \delta}) \notin I_{N+\log \delta}$ which means that $\sum_{k=1}^{N+\log \delta} 2^{-v_k} < \varepsilon_1$. For $k \in \{1, \ldots, N + \log \delta\}$ denote the middle quarter of $R_{2^{k+1}}$ by

$$\tilde{R}_{2^{k+1}} = \Xi_{2^{k+1}}\left(\left[-\frac{2^k}{2}, \frac{2^k}{2}\right]\right).$$

Let $\underline{\gamma}$ be a path chosen on the set of paths from $\tilde{R}_{2^k}$ to $\tilde{R}_{2^{k+1}}$ as follows. First, we uniformly choose a point $\underline{S}_k$ in $\tilde{R}_{2^k}$ and, independently, a point $\underline{A}_k$ in $\tilde{R}_{2^{k+1}}$. Given $\underline{S}_k$ and $\underline{A}_k$, we then fix $\gamma$ deterministically as with Lemma 3.3. Denote the probability involved in this selection procedure by $\tilde{P}$ and let $\hat{P} = P \otimes \tilde{P}$ denote the product measure of $P$ and $\tilde{P}$. Since $V(\gamma)^- \leq \sum_{e \in \gamma} X_e^-$ it follows that

$$\hat{E}[V(\gamma)^- \mid A(\mathbf{v})] \leq \hat{E}\left[\sum_{e \in T_k} X_e^- 1_{\{e \in \underline{\gamma}\}} \mid A(\mathbf{v})\right]$$
$$= \sum_{e \in T_k} E[X_e^- \mid A(\mathbf{v})] \tilde{P}(e \in \underline{\gamma}).$$

Furthermore, $\tilde{P}(e \in \underline{\gamma}) \leq c2^{-k}$ for some universal $c$, since $\underline{\gamma}$ was chosen just as in Lemma 3.3 and, on $A(\mathbf{v})$ [recall $A(\mathbf{v}) \subset B_k(v_k)$]

$$\sum_{e \in T_k} X_e \leq 2^{N-v_k+k+1}.$$

We thus have

$$\hat{E}[V(\gamma)^- \mid A(\mathbf{v})] \leq c2^{N-v_k+1}$$

and the Markov inequality gives

$$\hat{P}(V(\gamma)^- \geq 100c2^{N-v_k+1} \mid A(\mathbf{v})) \leq \frac{1}{100 2^{N-v_k+1}} E[V(\gamma)^- \mid A(\mathbf{v})]$$

$$\leq \frac{1}{100}.$$

In other words,

$$\hat{P}(V(\gamma)^- \geq 100c2^{N-v_k+1} \mid A(\mathbf{v})) = \hat{E}[\hat{E}[1_{\{V(\gamma)^- \geq -100 2^{N-v_k+1}\}} \mid \underline{S}_k] \mid A(\mathbf{v})]$$

$$\leq \tfrac{1}{100}$$

and, again using the Markov inequality, we have

$$\hat{P}(\hat{E}[1_{\{V(\gamma)^- \geq 100c2^{N-v_k+1}\}} \mid \underline{S}_k] \geq \tfrac{1}{10} \mid A(\mathbf{v}))$$

$$\leq 10 \hat{P}(V(\gamma)^- \geq 100c2^{N-v_k+1} \mid A(\mathbf{v}))$$

$$\leq \tfrac{1}{10}.$$

As $\underline{S}_k$ is uniformly chosen on $\tilde{R}_{2^k}$, on $A(\mathbf{v})$, the proportion of points $\underline{s}_k \in \tilde{R}_{2^k}$ verifying the inequality

$$\hat{P}(V(\gamma)^- \geq 100c2^{N+1-v_k} \mid \underline{S}_k = \underline{s}_k) \geq \tfrac{1}{10}$$

has to be less than $\tfrac{1}{10}$, that is, to say there exists a set, say $\mathcal{S}_{2^k} \subset \tilde{R}_{2^k}$, such that $\frac{|\mathcal{S}_{2^k}|}{|\tilde{R}_{2^k}|} \geq \tfrac{9}{10}$ and on $A(\mathbf{v})$, $\hat{P}(V(\gamma)^- \geq 100c2^{N-v_k+1} \mid \underline{S}_k \in \mathcal{S}_k) \leq \tfrac{1}{10}$. On the event $\{\underline{S}_k \in \mathcal{S}_{2^k}\}$, we can repeat the argument to obtain

$$\hat{P}(V(\gamma)^- \geq 100c2^{N+1-v_k} \mid A(\mathbf{v})) = \hat{E}[\hat{E}[1_{\{V(\gamma)^- \geq 100c2^{N+1-v_k}\}} \mid \underline{A}_k] \mid A(\mathbf{v})]$$

$$\leq \tfrac{1}{10}$$

and, once more using the Markov inequality, we obtain

$$\hat{P}(\hat{E}[1_{\{V(\gamma)^- \geq 100c2^{N-v_k+1}\}} \mid \underline{A}_k] \geq 1 \mid A(\mathbf{v}))$$

$$\leq \hat{P}(V(\gamma)^- \geq 100c2^{N-v_k+1} \mid \underline{S}_k \in \mathcal{S}_k \mid A(\mathbf{v}))$$

$$\leq \tfrac{1}{10}.$$

Again, as $\underline{A}_k$ was uniformly chosen independently of $\underline{S}_k$, this implies that there exists a set, say $\mathcal{A}_{2^k} \subset \tilde{R}_{2^{k+1}}$, such that $\frac{|\mathcal{A}_{2^k}|}{|\tilde{R}_{2^{k+1}}|} \geq \frac{9}{10}$ and

$$\hat{P}(V(\gamma)^- \geq 100c2^{n-v_k+1} \mid A_k \in \mathcal{A}_{2^k}, \underline{S}_k \in \mathcal{S}_{2^k}) < 1.$$

This implies that for every $\underline{s}_k \in \mathcal{S}_{2^k}$ and every $\underline{a}_k \in \mathcal{A}_{2^k}$, there exists a path $\gamma$ from $\underline{s}_k$ to $\underline{a}_k$ such that $V(\gamma)^- < 100c2^{N-v_k+1}$, which implies $V(\gamma) > -100c \times 2^{N-v_k+1}$.

Moreover, the construction of $\mathcal{S}_{2^k}$ and $\mathcal{A}_{2^k}$ is so that $|\mathcal{S}_{2^k} \cap \mathcal{A}_{2^{k-1}}| \geq \frac{7}{10}|\tilde{R}_{2^k}|$, and, therefore, if we take $\underline{s}_{k-1} \in \mathcal{S}_{2^{k-1}}, \underline{b}_k \in \mathcal{A}_{2^{k-1}} \cap \mathcal{S}_{2^k}$ and $\underline{a}_k \in \mathcal{A}_{2^k}$, there exists a path from $\underline{s}_{k-1}$ to $\underline{b}_k$ and a path from $\underline{b}_k$ to $\underline{a}_k$ whose concatenation is a path from $\underline{s}_{k-1}$ to $\underline{a}_k$ with value greater than $-100c2^{N-v_{k-1}+1} - 100c2^{N-v_k+1}$. So, for every $\underline{s}_{k-1} \in \mathcal{S}_{2^{k-1}}$ and every $\underline{a}_k \in \mathcal{A}_{2^k}$, there exists a path $\gamma$ from $\underline{s}_k$ to $\underline{a}_k$ such that $V(\gamma) > -100c2^{N+1}(2^{-v_{k-1}} + 2^{-v_k})$. Repeating the argument gives the existence of a path $\gamma$ from 0 to every point in $\mathcal{A}_{2^{N+\log \delta}}$ whose value satisfies

$$V(\gamma) > -c100 2^{N+1} \sum_{k=0}^{N+\log \delta} 2^{-v_k}.$$

Since $\sum_{i=1}^{N+\log \delta} 2^{-v_i} < \varepsilon_1$, we finally obtain that

$$V(\gamma) > -200c\varepsilon_1 2^N$$

and it suffices to choose $\varepsilon_1 = \frac{\varepsilon}{400c}$ to conclude the existence of a path from 0 to each point of $\mathcal{A}_{2^{N+\log \delta}} \subset \tilde{R}_{2^{N+\log \delta}}$ with the property that

$$V(\gamma) > -\frac{\varepsilon}{2} 2^N$$

and $|\mathcal{A}_{2^{N+\log \delta}}| \geq \frac{35}{50}|\tilde{R}_{2^{N+\log \delta}}|$. This gives us the existence of the desired paths. □

PROOF OF IMPLICATION PART OF THEOREM 1.1. Proposition 4.1 assures us that if event $A(\underline{v})$ occurs with $\sum_{k=0}^{N+\log \delta} 2^{-v_k}$, then for $\frac{7}{10}$ of the points, $\underline{x}$, in $\tilde{R}_{N+\log(\delta)}$ there exists a path, $\gamma^x$, from $(\underline{0},0)$ to $\underline{x}$ having value at least $-\frac{2^N \varepsilon}{100}$. Proposition 2.2, applied with $n = 2^N$, ensures that outside an event of probability $e^{-c(\varepsilon,\delta)2^{2N}}$ for $\frac{9}{10}$ of the points $\underline{x}$ in $\tilde{R}_{N+\log \delta}$ there are paths $\gamma'^{\underline{x}}$ from $\underline{x}$ to $H_{2^N}$ whose value is at least $(\mu - \varepsilon/10)(1-\delta)2^N$. Thus, for $\frac{9}{10}\frac{7}{10} - 1$ of the points $\underline{x}$ in $\tilde{R}_{N+\log \delta}$, there exist paths $\gamma^{\underline{x}}, \gamma'^{\underline{x}}$ with the requisite properties. But for such $\underline{x}$, the concatenation of these two paths gives a path from $(\underline{0},0)$ to $H_{2^N}$ of value at least $(\mu - \varepsilon/10)(1-\delta)2^N - \frac{2^N \varepsilon}{100}$. We conclude the result. Now a simple application of Proposition 2.1 taking $n = 2^N$, guarantees that except on an event of probability less than $e^{-C2^{2N}}$, there exists a set $B'_{2N+\log \delta} \subset \tilde{R}_{2N+\log \delta}$ so that $|B'_{2N+\log \delta}| \geq \frac{9}{20}|\tilde{R}_{2N+\log \delta}|$ and so that there exists from each point of $B'_{2N+\log \delta}$ a path $\underline{\gamma}_1$ leading to $R_{2^N}$ which has

$$V(\gamma_1) \geq \left(\mu - \frac{\varepsilon}{10}\right)(1-\delta)2^N.$$

Since $|B_{2N+\log \delta} \cap B'_{2N+\log \delta}| \geq \frac{15}{100}|\tilde{R}_{2n+\log \delta}|$, this shows that we can construct a path $\gamma \in \aleph_{2^N}$, by concatenation, which satisfies

$$V(\gamma) \geq \left(\mu - \frac{\varepsilon}{10}\right)(1-\delta)2^N - \frac{\varepsilon}{2}2^N$$
$$\geq (\mu - \varepsilon)2^N,$$

provided $\delta < \frac{2\varepsilon}{5\mu}$. $\square$

We have completed the proof of the implication: if (41), then

$$\varlimsup_{n \to \infty} \frac{1}{n^2} \log P(Z_n \leq (\mu - \varepsilon)n) < 0.$$

To complete the proof of Theorem 1.1, we now consider the case where condition (41) fails. We will show when (41) fails that, in fact, for all $\varepsilon > 0$,

(53) $$\lim_{n \to \infty} \frac{1}{n^2} \log P(Z_n \leq (\mu - \varepsilon)n) = 0.$$

As usual, we treat $n$ of the form $n = 2^N$. We first fix $\delta$ so that $2\delta M_0 < \varepsilon$. We will later also require that $\delta$ not exceed another constant. Now, we define for $N$ large and, in particular, $N \gg -\log \delta$,

$$\varepsilon_j \equiv 100\varepsilon \frac{1}{f(2^j)\sum_{k=-\log \delta}^{N} 1/(f(2^k))}, \qquad j = -\log \delta, \ldots, N.$$

Then $\sum_{j=-\log \delta}^{N} \varepsilon_j = 100\varepsilon$. We may suppose without loss of generality that $\varepsilon_j < 1$ for all $j = -\log \delta, \ldots, N$.

For reasons which will later become clear, we wish to consider just $j \in J$, where

$$J = \{k \in \{-\log \delta, \ldots, N\} : 2^k \varepsilon_k \geq M_0\}.$$

Note that

$$\sum_{\substack{j=-\log \delta \\ j \notin J}}^{N} \varepsilon_j < \sum_{j=-\log \delta}^{N} M_0 2^{-j}$$
$$< 2\delta M_0$$
$$< \varepsilon$$

and so
$$\sum_{j \in J} \varepsilon_j > 99\varepsilon.$$

Consider the event
$$V_J^N = \{\forall j \in J, X_e \leq -2^j \varepsilon_j, \forall e \in T_{N-j}\}.$$

Now by independence,
$$P(V_J^N) = \prod_{j \in J} \prod_{L(e)=N-j} P(X_e \leq -2^j \varepsilon_j)$$
$$\geq \prod_{j \in J} P(X_e \leq -2^j \varepsilon_j)^{2^{2(N-j+1)}}$$
(54)
$$= \exp\left\{-\sum_{j \in J} 4 \cdot 2^{2(N-j)}(2^j \varepsilon_j)^2 f(2^j \varepsilon_j)\right\}$$
$$= \exp\left\{-\sum_{j \in J} 4 \cdot 2^{2N} \varepsilon_j^2 f(2^j \varepsilon_j)\right\}$$
$$\geq \exp\left\{-\sum_{j \in J} 4 \cdot 2^{2N} \varepsilon_j^2 f(2^j)\right\}, \qquad f \text{ nondecreasing}$$
$$= \exp\left\{-4 \cdot 2^{2N} \sum_{j \in J} \varepsilon_j^2 f(2^j)\right\}.$$

Now from our definition of $\varepsilon_j$,
$$\varepsilon_j^2 f(2^j) = \frac{1}{f(2^j)(\sum_{k=-\log \delta}^{N} 1/(f(2^k)))^2},$$

so
$$\sum_{j \in J} \varepsilon_j^2 f(2^j) \leq \frac{1}{\sum_{k=-\log \delta}^{N} 1/(f(2^k))}.$$

By the positivity of $f$ we conclude
$$P(V_J^N) \geq \exp\left\{-2^{2N+2} \bigg/ \sum_{k=-\log \delta}^{N} \frac{1}{f(2^k)}\right\}.$$

Now by assumption $\sum_k \frac{1}{f(2^k)} = \infty$ and so as $N \to \infty$, $\sum_{k=-\log \delta}^{N} \frac{1}{f(2^k)}$ tends to $\infty$, which implies that
$$\lim_{N \to \infty} \frac{1}{2^{2N}} \log P(V_J^N) = 0.$$

Our result will hence be complete if we can show that

(55) $$P(Z_{2^N} \leq (\mu - \varepsilon)2^N \mid V_J^N) \geq 1/2.$$

We first consider the event

(56) $U_\delta^N = \{\forall \underline{\gamma}_1 \in \aleph_{\delta 2^N, 2^N}((-\delta 2^N, \delta 2^N)), V(\gamma_1) \leq (1-\delta)(\mu + \varepsilon)2^N\}.$

By (15), for $N$ large,

(57) $$P(U_\delta^N) \geq \tfrac{9}{10}.$$

Also, this event is independent of $V_J^N$. Secondly, we consider the behavior of paths in $\aleph_{\delta 2^N}$. The most important point is that for any $\underline{\gamma} \in \aleph_{\delta 2^N}$, and any $j \in J$, the path $\underline{\gamma}$ must contain $2^j$ edges of level $j$. Thus, on the event $V_J^N$, for any $\underline{\gamma} \in \aleph_{\delta 2^N}$ and $j \in J$

$$\sum_{e \in T_j, e \in \gamma} X_e \leq -2^j \varepsilon_j 2^{N-j}$$

$$= -2^N \varepsilon_j.$$

Consequently,

$$\sum_{j \in J} \sum_{e \in T_j, e \in \gamma} X_e \leq -2^N \sum_{j \in J} \varepsilon_j$$

$$\leq -99\varepsilon 2^N.$$

It remains to show that (with high probability) uniformly over possible $\underline{\gamma} \in \aleph_{\delta 2^N}$,

$$V(\gamma) \leq -10\varepsilon 2^N.$$

First, consider i.i.d. r.v.s $Y_e$ for $e \in T_j$ for some $j \leq N + \log \delta$ of distribution

$$P(Y_e \leq x) = P(X_e \leq x \mid X_e \geq -M_0).$$

We have that the $Y_e$ are stochastically greater than the $X_e$. We couple $\{X_e : e \in E\}$ and $\{Y_e : e \in E\}$ so that whenever $e \in T_j$ with $j \notin J, X_e \leq Y_e$ and $\sigma(\{Y_e : e \in T_j, j \in J\})$ is independent of $\sigma(\{X_e : e \in T_j, j \in J\})$. We have, for some constant $\mu'$ not depending on $\delta$ or $N$ (just on the distribution of $Y_e$) that as $N \to \infty$,

$$\frac{\sup_{\underline{\gamma}} \sum_{e \in \gamma} Y_e}{\delta 2^N} \xrightarrow{\text{pr}} \mu' < \infty,$$

where the supremum is taken over paths $\underline{\gamma} \in \aleph_{\delta 2^N}$. Therefore, the event

$$W \equiv \left\{ \sup_{\underline{\gamma}} \frac{\sum_{e \in \gamma} Y_e}{\delta 2^N} < (\mu' + 1) \right\}$$

satisfies (by independence) for $N$ large

(58) $$P(W \mid V_J^N) > 9/10.$$

Then on the event $W$, uniformly in $\gamma$,

$$\sum_{e \in \gamma} X_e = \sum_{e \in \gamma, e \in T_j, j \in J} X_e + \sum_{e \in \gamma, e \in T_j, j \notin J} X_e$$

$$\leq \sum_{e \in \gamma, e \in T_j, j \in J} X_e + \sum_{e \in \gamma, e \in T_j, j \notin J} Y_e$$

$$\leq \sum_{e \in \gamma, e \in T_j, j \in J} X_e + \sum_{e \in \gamma} Y_e - \sum_{e \in \gamma, e \in T_j, j \in J} Y_e.$$

On $V_J^N$, the first term is less than $-99\varepsilon 2^N$. On $W$, the second term is less than $(\mu' + 1)\delta 2^N$. Finally, since necessarily $Y_e \geq -M_0$ for all $e$, the third term satisfies

$$\sum_{e \in \gamma, e \in T_j, j \in J} Y_e \leq M_0 \delta 2^N.$$

Hence, on $V_J^N \cap W$ for all $\underline{\gamma} \in \aleph_{\delta 2^N}$,

$$V(\gamma) \leq -99\varepsilon 2^N + 2^N \delta(\mu' + 1 + M_0)$$
$$\leq -10\varepsilon 2^N$$

provided $\delta$ was chosen sufficiently small.

Thus, on $V_J^N \cap W \cap U_\delta^N$, we have for any $\underline{\gamma} \in \aleph_{2^N}$

$$V(\gamma) = \sum_{e \in \gamma, e \in T_j, j \leq N + \log \delta} X_e + \sum_{e \in \gamma, e \in T_j, j \geq N + \log \delta} X_e$$

(59) $$\leq -10\varepsilon 2^N + (1 - \delta)(\mu + \varepsilon)2^N$$
$$\leq (\mu - \varepsilon)2^N.$$

Combining (57), (58) and (59) gives (55) and we are done.

**5. Model 2.** Finally, we give a sketch for $d = 2$ of how the proofs are adapted to the situation in Model 2. There are two minor differences between the two models, the first is that the underlying graph structures are slightly different. The second is that in Model 1 the trouble arises from large negative values in the field, whereas in Model 2, the trouble arises from large positive values of the field. The analysis, nonetheless, proceeds in a similar manner. For the purposes of establishing the upper bound claimed in Theorem 1.3,

(60) $$\varlimsup_{n \to \infty} \frac{1}{n^2} \log P(a_n > (\nu + \varepsilon)n) < 0,$$

we can consider the infimum over a smaller set of paths. As before, we will consider, purely for notational convenience, $a_{0,n}$ for $n$ of the form $2^N$. Write the typical path in coordinates as $\gamma = (\gamma_1, \gamma_2)$. Next, take $\delta > 0$ small and consider the set of paths defined by

(61)
$$\begin{aligned}
\Psi_{\delta,N} = \{\gamma \in \Psi_{2^N} : &\gamma(0) = (0,0), \gamma(l(\gamma)) = (2^N, 0), \\
&\gamma_1(j+1) \geq \gamma_1(j), \text{ whenever } \gamma_1(j) \in [0, \delta 2^N) \cup [(1-\delta)2^N, 2^N] \\
&\gamma_1(j) \geq |\gamma_2(j)|, \text{ for } \gamma_1(j) \in [0, \delta 2^N), \\
&2^N - \gamma_1(j) \geq |\gamma_2(j)|, \text{ for } \gamma_1(j) \in [(1-\delta)2^N, 2^N], \\
&\gamma_2(j) \in [-\delta 2^N, \delta 2^N], \text{ for } \delta 2^N \leq \gamma_1(j) \leq (1-\delta)2^N\}.
\end{aligned}$$

Next, set

(62) $$a_{\delta,2^N} = \inf_{\gamma \in \Psi_{\delta,N}} \sum_{e \in \gamma} t_e.$$

Then, obviously,

(63) $$P(a_{\delta,2^N} \geq (\nu + \varepsilon)2^N) \geq P(a_{2^N} \geq (\nu + \varepsilon)2^N)$$

and so we can establish (60) by proving

(64) $$\varlimsup_{2^N \to \infty} \frac{1}{2^{2N}} \log P(a_{\delta,2^N} > (\nu + \varepsilon)2^N) < 0.$$

While in the region $\{x \in \mathbb{Z}^2 : 0 \leq x_1 \leq \delta 2^N\} \cup \{x \in \mathbb{Z}^2 : (1-\delta)2^N \leq x_1 \leq 2^N\}$, paths in $\Psi_{\delta,2^N}$ behave very nearly the same as in Model 1. Indeed, the graphical and path structures are quite similar. One now considers the "shells" $T_k$ of edges with distance from the starting or terminal point between $2^k$ and $2^{k+1}$, for $k = 0, 1, \ldots, N + \log \delta$. More precisely,

$$T_k = \{e = \{x,y\} : 2^k \leq |x_1| < 2^{k+1}\} \cup \{e = \{x,y\} : 2^k \leq |2^N - y_1| < 2^{k+1}\}.$$

Nearly the same arguments used in the case of Model 1 can be used here. The only change is that one deals directly with $\sum_{e \in T_k} t_e$ instead of $V_k^-$ as defined at (29). In the present case, one considers

(65) $$\tilde{Z}_k = \sum_{e \in T_k} t_e.$$

The main concern is still with the event that the sum $\sum_{e \in T_k} t_e$ may be large. Thus, similarly to (66), we define

(66) $$\tilde{B}_k(v) = \{2^{-k} Z_k \in [2^{N-v}, 2^{N-v+1})\}$$

and corresponding versions $\tilde{A}(\mathbf{v})$ and $\tilde{A}^I(\mathbf{v})$ as in (43) and (44). Then an appropriately modified version of Proposition 4.1 follows using the estimates derived as in Lemma 4.1 and Lemma 4.2. The random selection of path procedure leads to the same probability estimate for an edge to be on a randomly selected path as in the proof of part (2) of Proposition 4.1. We now substitute Proposition 2.2 for Proposition 2.1. Then as in the proof of the implication part of Theorem 1.1, concatenating paths from 0 to the face $\{x \in \mathbb{Z}^2 : x_1 = \delta 2^N, |x_2| \leq \delta 2^N\}$ with value at most $\frac{\varepsilon}{100} 2^N$ from the face $\{x \in \mathbb{Z}^2 : x_1 = \delta 2^N, |x_2| \leq \delta 2^N\}$ to the face $\{x \in \mathbb{Z}^2 : x_1 = (1-\delta) 2^N, |x_2| \leq \delta 2^N\}$ with value at most $(\nu + \frac{\varepsilon}{5})(1 - 2\delta) 2^N$ and from the face $\{x \in \mathbb{Z}^2 : x_1 = (1-\delta) 2^N, |x_2| \leq \delta 2^N\}$ to $(2^N, 0, \ldots, 0)$ with value at most $\frac{\varepsilon}{100} 2^N$ all with probability exceeding $1 - e^{-c 2^{2N}}$ gives that

$$P(a_{2^N} \leq (\nu + \varepsilon) 2^N) \geq 1 - e^{-c 2^{2N}}.$$

We note that in choosing the oriented lattice structure for our paths at the start and end, we are giving up a lot. However, it must be borne in mind that the objective is to find some $\delta$ so that there will, outside an event of very small probability, exist many paths from $(0, 0)$ to the interval $\{2^{N\delta}\} \times [-\delta 2^N, \delta 2^N]$ whose values are not larger than $\varepsilon 2^N / 2$. We do not try to find an optimal $\delta$. For the "only if" claim of Theorem 1.3, namely,

(67) $$\lim_{N \to \infty} \frac{1}{2^{Nd}} \log P(a_{2^N} > (\nu + \varepsilon) 2^N) = 0$$

holds when

(68) $$\sum_{k=1}^{\infty} \frac{1}{f(2^k)^{1/(d-1)}} = \infty$$

the argument follows the same lines. There are slight differences. For example, again restricting to $d = 2$ for convenience, in the present context, define the level of an edge $e$ to be $L(e) = k$ when $e = \{x, y\}, x, y \in \mathbb{Z}^2$ and $|x|_1 \wedge |y|_1 = i \in [2^k, 2^{k+1})$, where $|x|_1 = |x_1| + |x_2|$. Then put

$$W_J^N = \{\forall j \in J, t_e \geq 2^j \varepsilon_j, \forall e \text{ with } L(e) = N - j\}.$$

Then exactly as in (54) one has

(69) $$P(W_J^N) \geq \exp\left\{-16 \cdot 2^{2N} \sum_{j \in J} \varepsilon_j^2 f(2^j)\right\}.$$

This yields

$$P(W_J^N) \geq \exp\left\{-2^{2N+4} \Big/ \sum_{k=-\log \delta}^{N} \frac{1}{f(2^k)}\right\}.$$

A result analogous to (58) follows from (14). More precisely, define for appropriately small $\delta$,

(70) $$\Theta_{\delta, N} = \{\gamma \in \Theta : \gamma(0) \in \partial[-\delta 2^N, \delta 2^N]^2, \gamma(l(\gamma)) = (2^N, 0)\}$$

and set

(71) $$\tilde{U}_\delta^N = \{\forall \gamma \in \Theta_{\delta,2^N}, V(\gamma) \geq (1-\delta)(\nu-\varepsilon)2^N\}.$$

We claim that for $N$ large,

(72) $$P(\tilde{U}_\delta^N) \geq \tfrac{9}{10}.$$

To see this, first observe that for $N$ large,

$$P(\tilde{U}_0^N) \geq \tfrac{99}{100}.$$

Consider the event, $E_{\delta 2^N}$, that for some random point $p \in \partial[-\delta 2^N, \delta 2^N]^2$ we have a path $\gamma_p \in \Theta_{\delta,N}$ such that $\gamma_p(0) = p$ and $V(\gamma) < (1-\delta)(\nu-\varepsilon)2^N$, which is less than $(\nu - \varepsilon/2)2^N$ if $\delta$ is chosen sufficiently small. We may assume that $\gamma_p(k) \notin (-\delta 2^N, \delta 2^N)^2$ $\forall k > 0$, otherwise we could wait until the last exit from $[-\delta 2^N, \delta 2^N]^2$ and proceed from there on a path with a smaller value. Thus, $V(\gamma_p)$ is independent of $\{t_e : e \in (-\delta 2^N, \delta 2^N)^2\}$. Consequently, we can take the shortest path with no more than $2\delta 2^N$ edges from 0 to $p$. Call this path $\hat{\gamma}_p$. By the law of large numbers, we can select $N$ large enough so that $P(|V(\hat{\gamma}_p) - E[t_e]| \leq \varepsilon \delta 2^N) \geq 99/100$. Concatenating $\hat{\gamma}_p$ and $\gamma_p$ gives a path $\gamma$ with $V(\gamma) \leq (1-\delta)(\nu-\varepsilon)2^N + E[t_e] + \varepsilon \delta 2^N$. The remainder of the argument follows just as in the case of Model 1.

M. CRANSTON
DEPARTMENT OF MATHEMATICS
UNIVERSITY OF CALIFORNIA, IRVINE
IRVINE, CALIFORNIA 92612
USA
E-MAIL: mcransto@math.uci.edu

D. GAUTHIER
T. S. MOUNTFORD
DMA
ÉCOLE POLYTECHNIQUE FÉDÉRALE DE LAUSANNE
CH-1015
SWITZERLAND
E-MAIL: damien.gauthier@epfl.ch
        thomas.mountford@epfl.ch